\newcommand{\siamyesno}[2]{#2}   
\siamyesno{
\documentclass[onetabnum,onefignum,nohypdvips,final]{siamart171218}
\usepackage{amssymb,amsmath,epsfig,verbatim}
\newtheorem{remark}[theorem]{Remark}
\newtheorem{ass}[theorem]{Assumption}
}{
\documentclass[11pt,a4paper]{article}
\usepackage{amssymb,amsmath,amsthm,epsfig,verbatim,xcolor}
\newtheorem{theorem}{\sc Theorem.}[section]
\newtheorem{lemma}[theorem]{\sc Lemma.}
\newtheorem{remark}[theorem]{\sc Remark.}
\newtheorem{ass}[theorem]{\sc Assumption.}

\newenvironment{AMS}%
{{\upshape\bfseries AMS subject classifications. }\ignorespaces}{}
\newenvironment{keywords}{{\upshape\bfseries Key words. }\ignorespaces}{}

}
\usepackage{mhchem} 
\usepackage{ifpdf}  
\ifpdf
  \DeclareGraphicsExtensions{.eps,.pdf,.png,.jpg}
\else
  \DeclareGraphicsExtensions{.eps,.ps}
\fi

\newcommand{\bRgeq}{{\mathbb R}_{\geq 0}}
\newcommand{\RZ}{{\mathbb R} \slash {\mathbb Z}}
\newcommand{\bR}{{\mathbb R}}

\newcommand{\matpartnpar}[1]{\partial_{#1}^{\text{\tiny$\square$}}}
\newcommand{\matpartn}{\matpartnpar{t}}
\newcommand{\erf}{\operatorname{erf}}

\newcommand{\drho}{\;{\rm d}\rho}
\newcommand{\du}{\;{\rm d}u}

\newcommand{\dd}[1]{\frac{\rm d}{{\rm d}#1}}
\newcommand{\ddt}{\dd{t}}

\def\epsilon{\varepsilon} 
\def\hat{\widehat}

\newcommand{\errorXtxtc}{\max_{m=1,\ldots,M-1}\| \dot x(t_m) - \tfrac{ x^{m+1} -  x^{m-1}}{2\Delta t} \|_{0,h}}
\siamyesno{}{
\textwidth 455pt \oddsidemargin 0pt \evensidemargin 0pt \headsep
0pt \headheight 0pt \textheight 655pt \parskip 10pt \parindent 0pt
}

\begin{document}
\title{
Discrete hyperbolic curvature flow in the plane
}
\author{Klaus Deckelnick\footnotemark[2]\ \and 
        Robert N\"urnberg\footnotemark[3]}

\renewcommand{\thefootnote}{\fnsymbol{footnote}}
\footnotetext[2]{Institut f\"ur Analysis und Numerik,
Otto-von-Guericke-Universit\"at Magdeburg, 39106 Magdeburg, Germany \\
{\tt klaus.deckelnick@ovgu.de}}
\footnotetext[3]{Dipartimento di Mathematica, Universit\`a di Trento,
38123 Trento, Italy \\ {\tt robert.nurnberg@unitn.it}}

\date{}

\maketitle

\begin{abstract}
Hyperbolic curvature flow is a geometric evolution equation that in the plane
can be viewed as the natural hyperbolic analogue of curve shortening flow.
It was proposed by Gurtin and Podio-Guidugli (1991) to model certain wave
phenomena in solid-liquid interfaces.
We introduce a semidiscrete finite difference method for the 
approximation of hyperbolic curvature flow and prove error bounds for 
natural discrete norms. We also
present numerical simulations, including the onset of singularities starting
from smooth strictly convex initial data.
\end{abstract} 


\renewcommand{\thefootnote}{\arabic{footnote}}

\siamyesno{
\pagestyle{myheadings}
\thispagestyle{plain}
\markboth{K. DECKELNICK AND R. N\"URNBERG}
{DISCRETE HYPERBOLIC CURVATURE FLOW IN THE PLANE}
}{}

\setcounter{equation}{0}
\section{Introduction} 

The analytical and numerical study of parabolic geometric evolution equations,
such as mean curvature flow, surface diffusion and Willmore flow, to name a
few, has received considerable attention in the literature over the last 
few decades, see e.g.\
\cite{Mullins57,Huisken84,GageH86,Willmore93,TaylorC94,ElliottG97a,Ecker04,%
DeckelnickDE05,Giga06,Mantegazza11,bgnreview}. 
On the other hand, hyperbolic evolution laws for moving interfaces have been
studied far less. In this paper, we are going to investigate the numerical
approximation of the hyperbolic geometric evolution equation
\begin{equation} \label{eq:GurtinGamma}
\alpha \matpartn \mathcal{V}_{\Gamma} + \beta \mathcal{V}_{\Gamma}
= \varkappa_{\Gamma} \quad \text{on }\ \Gamma(t),
\end{equation}
for a family of closed curves $(\Gamma(t))_{t\in[0,T]}$ in $\bR^2$. Here 
$\mathcal{V}_\Gamma$ denotes the velocity of $(\Gamma(t))_{t\in[0,T]}$
in the direction of the normal $\nu_\Gamma$,
$\matpartn$ is the normal time derivative on $(\Gamma(t))_{t\in[0,T]}$,
and $\varkappa_\Gamma$ denotes the curvature of $\Gamma(t)$. 
Our sign convention is such that the unit circle with outward normal has 
curvature $\varkappa_\Gamma = -1$.
The flow \eqref{eq:GurtinGamma} corresponds to the evolution law proposed in
\cite[(1.2)]{GurtinP91}, in the case of an isotropic surface energy and in 
the absence of external forcings, where it was suggested as a model for the
evolution of melting-freezing waves at the solid-liquid interface of crystals
such as \ce{^{4}He} helium. Here the parameters
$\alpha \in \bRgeq$ and $\beta \in \bRgeq$ play the role of an effective
density and a kinetic coefficient, respectively. 
In the special case
$\alpha=1$ and $\beta=0$ we obtain the hyperbolic geometric evolution law
\begin{equation} \label{eq:Gurtinbeta0}
\matpartn \mathcal{V}_{\Gamma} = \varkappa_{\Gamma} \quad\text{on }\ \Gamma(t),
\end{equation}
while the choices 
$\alpha=0$ and $\beta=1$ yield  the well-known (mean) curvature flow, or
curve shortening flow. However, since in this work we are interested in the 
hyperbolic case, we shall from now on set $\alpha=1$ for simplicity.
We remark that in order to close the geometric evolution equation
\eqref{eq:GurtinGamma}, the initial conditions
\[
\Gamma(0) = \Gamma_0 \quad \text{and}\quad
\mathcal{V}_\Gamma\!\mid_{t=0} = \mathcal V_{\Gamma,0}
\]
need to be prescribed, where $\Gamma_0$ defines the initial curve and
$\mathcal V_{\Gamma,0} : \Gamma_0 \to \bR$ gives an initial normal velocity.

Let us consider a parametric description of the evolving curves, i.e.\ 
$\Gamma(t)=  x(I,t)$ for some mapping $ x:I \times [0,T] \to \bR^2$,
where $I=\RZ$ is the periodic interval $[0,1]$. We denote by 
\begin{equation} \label{eq:nu}
\tau =  \frac{x_\rho}{| x_\rho |}, \quad
\nu = \tau^\perp =  \frac{x_\rho^\perp}{| x_\rho|}  \quad  \text{and} \quad \varkappa \nu = \frac{\tau_\rho}{| x_\rho |} =
\frac{1}{| x_\rho |}\bigl(\frac{x_\rho}{| x_\rho |} \bigr)_\rho ,
\end{equation}
the unit tangent, the unit normal and the curvature vector, respectively, 
so that e.g.\ $\nu = \nu_\Gamma \circ x$ and 
$\varkappa=\varkappa_\Gamma \circ x$. Here and throughout
$\cdot^\perp$ denotes the anti-clockwise rotation through 
$\frac{\pi}{2}$.
We shall show in
Lemma~\ref{lem:normalflow} below that if $x$ is a solution of the system
\begin{subequations} \label{eq:hcsf}
\begin{alignat}{2}
x_{tt} + \beta  x_t & =  \frac{1}{| x_\rho |}
\bigl(\frac{x_\rho}{| x_\rho |} \bigr)_\rho - ( x_t \cdot \tau_t) \tau 
&&\quad \text{in } I \times (0,T], \label{eq:abxttxss} \\
 x(\cdot, 0) & =  x_0,\ 
 x_t(\cdot, 0) = \mathcal{V}_0\nu(\cdot,0) &&\quad \text{in } I, \label{eq:init}
\end{alignat}
\end{subequations}
then the curves $(\Gamma(t))_{t \in [0,T]}$  evolve according to 
\eqref{eq:GurtinGamma} with $\alpha=1$. In the above, 
$x_0 : I \to \bR^2$ is a parameterization of the given initial 
curve $\Gamma_0$
and $\mathcal{V}_0 =\mathcal{V}_{\Gamma,0} \circ x_0$ 
is induced by the given initial normal velocity $\mathcal{V}_{\Gamma,0}$.
The introduction of the second term
on the right hand side of \eqref{eq:abxttxss} has the effect that
the parameterization $x$ is normal,
i.e.\ it satisfies $x_t \cdot \tau =0$, see also Lemma~\ref{lem:normalflow}.
The system \eqref{eq:hcsf} in the case $\beta=0$ 
has been studied in \cite{KongLW09,KongW09}, see also \cite{HeHX17}. 
In particular, it is shown in \cite{KongLW09}
that if $\Gamma(0)$ is strictly convex,  
and if the initial velocity $\mathcal V_0 \nu(\cdot,0)$ does not point
outwards anywhere on $\Gamma(0)$,
then the solution to \eqref{eq:hcsf}
exists on a finite time interval $[0,T_{\max})$ and the curves $\Gamma(t)$ remain strictly convex. 
Furthermore, as $t \rightarrow T_{\max}$, 
$\Gamma(t)$ either shrinks to a point or converges to a convex curve with discontinuous curvature.

One may wonder whether it is possible to replace \eqref{eq:abxttxss} by
the simpler hyperbolic equation
\begin{equation} \label{eq:xttxss0}
 x_{tt} =  \frac{1}{| x_\rho |}
\bigl(\frac{x_\rho}{| x_\rho |} \bigr)_\rho \quad \text{in } I \times (0,T],
\end{equation}
which has been considered in e.g.\ \cite{HeKL09} after having been 
proposed by Yau in \cite[p.~242]{Yau00}.
However, in contrast to \eqref{eq:hcsf}, it is not clear whether solutions
to \eqref{eq:xttxss0} with the initial conditions \eqref{eq:init} parameterize
solutions to the flow \eqref{eq:Gurtinbeta0}. 
In fact, numerical evidence in Section~\ref{sec:nr54}, below, suggests
that solutions to \eqref{eq:hcsf} and \eqref{eq:xttxss0}, \eqref{eq:init}
parameterize different curve evolutions.
 
An alternative hyperbolic geometric evolution equation, 
that is similar to \eqref{eq:hcsf}, and which has been considered in 
\cite{LeFlochS08}, is described by
\begin{equation} \label{eq:LeFloch0}
 x_{tt} = \tfrac12 ( | x_t|^2 + 1)  \frac{1}{| x_\rho |}
\bigl(\frac{x_\rho}{| x_\rho |} \bigr)_\rho 
- ( x_t \cdot \tau_t) \tau \quad\! \text{in } I \times (0,T],
\quad
x(\cdot, 0) = x_0,\ x_t(\cdot, 0) = \mathcal{V}_0\nu(\cdot,0)
\quad\! \text{in } I.
\end{equation}
It can be shown that solutions to \eqref{eq:LeFloch0} also represent normal
parameterizations of curves. 
An interesting aspect of \eqref{eq:LeFloch0} in terms of the analysis
is that its solutions satisfy the energy conservation
\begin{equation*} 
\tfrac12 \ddt \int_I (| x_t|^2 + 1) | x_\rho| \drho = 0.
\end{equation*}
In contrast, for the flow \eqref{eq:hcsf} a conditional decay property can be
shown for the energy $\tfrac12\int_I (| x_t|^2 + 2) | x_\rho| \drho$, see
Remark~\ref{rem:normalflow} below, something that we will utilize for the
numerical analysis presented in this paper.
Let us finally mention that geometric second order hyperbolic PDEs have 
recently been used in \cite{DongHZ21} for applications in image processing. 

As regards the numerical approximation of hyperbolic geometric evolution
equations in the literature, we are only aware of the works
\cite{RotsteinBN99} and \cite{GinderS16}. In the former an algorithm for the
evolution of polygonal curves under crystalline hyperbolic curvature flow is
presented, which corresponds to \eqref{eq:GurtinGamma} for a crystalline,
anisotropic surface energy. 
On the other hand, in \cite{GinderS16} a level-set approach,
which is based on a threshold algorithm of BMO type, is used
for the numerical solution of \eqref{eq:xttxss0}. 

In this paper we will present a finite difference approximation of
\eqref{eq:hcsf} and prove an error bound for it. 
To the best of our knowledge this is the 
first result on the numerical analysis for a hyperbolic geometric evolution
equation in the literature. 

The remainder of the paper is organized as follows. In Section~\ref{sec:mf}
we show that curves $\Gamma(t)$ that are parameterized by solutions of 
\eqref{eq:hcsf} evolve according to \eqref{eq:GurtinGamma}.
We also derive several properties of these solutions.  
In Section~\ref{sec:fd} we
introduce our semidiscrete finite difference approximation and 
state our main result, Theorem~\ref{thm:main}. Its proof is presented in
Section~\ref{sec:ee}.
Finally, in Section~\ref{sec:nr} we suggest a fully discrete scheme
and present several numerical simulations for it,
including a convergence experiment
and simulations that lead to nonvanishing singularities in finite time.

\setcounter{equation}{0}
\section{Mathematical formulation} \label{sec:mf}

Consider a family $(\Gamma(t))_{t\in[0,T]}$ of evolving curves that are given by $\Gamma(t) = x(I,t)$, where
$ x : I \times [0,T] \to \bR^2$ satisfies $| x_\rho| > 0$ in $I \times [0,T]$.   Then the unit normal on $\Gamma$,
the curvature of $\Gamma$, the normal velocity of $\Gamma$ as well as the
normal time derivative on $\Gamma$ are defined by the following 
identities in $I$, see e.g.\ \cite{bgnreview}:
\begin{align} \label{eq:ItoGamma}
\nu_\Gamma \circ  x = \nu, \quad
\varkappa_\Gamma \circ  x = \varkappa, \quad
\mathcal{V}_\Gamma \circ  x =  x_t \cdot \nu,\quad 
(\matpartn f) \circ  x = 
(f \circ  x)_t - (f \circ  x)_s  x_t \cdot \tau,
\end{align}
where $\partial_s = | x_\rho |^{-1} \partial_\rho$ denotes differentiation with respect to arclength $s$.
We stress that the definitions of the above
quantities are independent of the chosen parameterization. 
The following lemma
establishes the connection to the evolution law
\eqref{eq:GurtinGamma} and derives additional properties of $x$ that will be useful in the subsequent analysis.

\begin{lemma} \label{lem:normalflow}
Suppose that $ x : I \times [0,T] \to \bR^2$ is a solution of \eqref{eq:hcsf}. Then the curves $(\Gamma(t))_{t \in [0,T]}$ 
with $\Gamma(t)=x(I,t)$ evolve according to \eqref{eq:GurtinGamma}. 
Furthermore, $x$ is a normal parameterization, i.e.
\begin{equation} \label{eq:xnormal}
x_t \cdot \tau =0 \quad \mbox{ in } I \times [0,T]
\end{equation}
and satisfies 
\begin{equation} \label{eq:leequation}
\partial_t | x_\rho |   = - |x_\rho | \, x_t \cdot x_{tt} - \beta  | x_\rho | \, | x_t |^2 \qquad \mbox{ in } I \times [0,T].
\end{equation}
\end{lemma}
\begin{proof} 
Using \eqref{eq:abxttxss} and \eqref{eq:nu} we deduce that
\begin{equation*} 
( x_t \cdot \tau)_t 
= x_{tt} \cdot \tau +   x_t \cdot \tau_t 
= ( \varkappa \nu -  ( x_t \cdot \tau_t) \tau - \beta  x_t) 
\cdot\tau +   x_t \cdot \tau_t 
= -\beta  x_t\cdot\tau.
\end{equation*}
In view of \eqref{eq:init} we have  $( x_t \cdot \tau)\!\mid_{t=0} = 0$ which implies \eqref{eq:xnormal}. With 
the help of  \eqref{eq:ItoGamma} and  \eqref{eq:xnormal} we now deduce
\begin{align*}
(\matpartn \mathcal{V}_{\Gamma}) \circ x + \beta \mathcal{V}_{\Gamma} \circ x 
& = [( x_t \cdot \nu)_t - ( x_t \cdot \nu)_s  x_t \cdot \tau] 
+ \beta  x_t \cdot \nu \\ & 
= x_{tt} \cdot \nu + x_t \cdot \nu_t + \beta x_t \cdot \nu  
= x_{tt} \cdot \nu + \beta x_t \cdot \nu \\ &
= \varkappa \nu \cdot \nu = \varkappa = \varkappa_\Gamma \circ x 
\qquad \mbox{ in } I \times [0,T],
\end{align*}
where we used that $0=\tfrac12 (| \nu |^2)_t= \nu_t \cdot \nu$, \eqref{eq:nu} and \eqref{eq:abxttxss}. Thus \eqref{eq:GurtinGamma} holds on
$\Gamma(t)$.  Finally, recalling again \eqref{eq:xnormal} and \eqref{eq:abxttxss}, we obtain
\begin{equation} \label{eq:leequation1}
\partial_t | x_\rho |  = x_{t \rho} \cdot \tau = - x_t \cdot \tau_\rho =  - |x_\rho | \, x_t \cdot x_{tt} - \beta   | x_\rho | \, | x_t |^2 \qquad \mbox{ in } I \times [0,T],
\end{equation}
which proves \eqref{eq:leequation}.
\end{proof}

\begin{remark} \label{rem:normalflow}
Using \eqref{eq:leequation1} and \eqref{eq:xnormal} we derive the following 
energy law
\begin{align} \label{eq:lemdte}
\tfrac12 \ddt \int_I (|x_t|^2 + 2) | x_\rho| \drho &
= \tfrac12 \int_I | x_t |^2  \partial_t | x_\rho | \drho + \int_I x_t \cdot x_{tt}  | x_\rho |+ \partial_t | x_\rho |
\drho \nonumber \\ &
= - \tfrac12 \int_I | x_t |^2 x_t \cdot \tau_\rho \drho  - \beta \int_I | x_t |^2 | x_\rho | \drho 
\nonumber \\ &
= - \tfrac12 \int_I (x_t \cdot \nu)^3 \varkappa | x_\rho | \drho - \beta \int_I (x_t \cdot \nu)^2 | x_\rho | \drho,
\end{align}
which corresponds to \cite[(4.6)]{GurtinP91} in the absence of external forces.
An adaptation of this relation to the error between continuous and discrete 
solution will be at the heart of our error analysis.
\end{remark}

For the remainder of the paper we make the following regularity assumptions concerning the solution $x$.

\begin{ass} \label{ass:x}
$x : I \times [0,T] \to \bR^2$ is a solution of \eqref{eq:hcsf} such that
 $\partial_t^i \partial_\rho^j x$ exist and are continuous on $I \times [0,T]$ for all $i,j \in \mathbb N \cup \{0\}$ with $2i+j \leq 4$. 
Furthermore, $|x_\rho| > 0$ in $I \times [0,T]$.
\end{ass}
 
Assumption~\ref{ass:x} implies in particular
that there exist constants $0< c_0 \leq C_0$ such that
\begin{equation} \label{eq:lereg}
c_0 \leq | x_\rho | \leq C_0 \quad \mbox{ in } I \times [0,T],\qquad 
\max_{I \times [0,T]} \left(
|\tau_\rho| + | x_t | + | x_{t \rho} | \right) \leq C_0 .
\end{equation}

\setcounter{equation}{0}
\section{Finite difference discretization} \label{sec:fd}

We shall employ a finite difference scheme in order to discretize 
\eqref{eq:hcsf} in space. 
To do so, let us introduce the set of grid 
points $\mathcal G^h:= \{\rho_1,\ldots, \rho_J \} \subset I$, where 
$\rho_j = jh$, $j=0,\ldots, J$, and $h = \frac 1J$ for $J\geq2$. 
In order to account for our periodic setting we always identify 
$\rho_0$ with $\rho_J$.
For a grid function $v:\mathcal G^h \to \bR^2$ we write 
$v_j:= v(\rho_j)$, $j=1,\ldots,J$, and in addition set
$v_0 = v_J$ and $v_{J+1} = v_1$ in view of the periodicity of $I$.
We associate with $v$ the backward difference quotient:
\begin{equation} \label{eq:delta-}
\delta v_j:= \frac{v_j - v_{j-1}}{h},
\quad j=1,\ldots,J
\end{equation}
and introduce the following discrete norms
\begin{equation} \label{eq:discnorms}
\| v \|_{0,h}: = \Bigl( h \sum_{j=1}^J | v_j |^2 \Bigr)^{\frac{1}{2}}, \quad 
\| v \|_{1,h}:= \Bigl( h \sum_{j=1}^J \bigl( | v_j |^2 + | \delta v_j |^2 \bigr) \Bigr)^{\frac{1}{2}}.
\end{equation}

Let $x^h: \mathcal G^h \to \bR^2$ be a grid function that will play the role of
a discrete parameterization of a curve. Then on $I_j=[\rho_{j-1},\rho_j]$, 
the associated discrete length element $q^h_j$ and the discrete tangent 
$\tau^h_j$ are given by 
\begin{equation*} 
q^h_j = | \delta x^h_j |, \quad \tau^h_{j} = \frac{1}{q^h_j} \delta x^h_j, 
\quad j=1,\ldots,J.
\end{equation*}
It will be convenient to also introduce the averaged vertex tangent 
$\theta^h_j$ via
\begin{equation} \label{eq:fdthetah}
\theta^h_j = 
\frac{\tau^h_{j} + \tau^h_{j+1}}{|\tau^h_{j} + \tau^h_{j+1}|}, \quad \mbox{ provided that } \tau^h_j + \tau^h_{j+1} \neq 0, \quad j=1,\ldots,J.
\end{equation}
Clearly, 
\begin{equation} \label{eq:fdthetaorthogonal}
(\tau^h_{j+1} - \tau^h_{j}) \cdot \theta^h_j
 = (\tau^h_{j+1} - \tau^h_{j}) \cdot \frac{\tau^h_{j} + \tau^h_{j+1}}{|\tau^h_{j} + \tau^h_{j+1}|} = \frac1{|\tau^h_{j} + \tau^h_{j+1}|}
( |\tau^h_{j+1}|^2 - |\tau^h_{j}|^2) = 0.
\end{equation}

\begin{lemma} \label{lem:prop1}
Let $x \in C^4(I;\bR^2)$ such that $c_0 \leq | x_\rho | \leq C_0$ in $I$ 
and set $\tau = \frac{x_\rho}{|x_\rho|}$ as well as
\begin{displaymath}
x_j=x(\rho_j), \; q_j = | \delta x_j|, \; \mbox{ and } \tau_{j} = \frac{1}{q_{j}} \delta x_{j}, \; j=1,\ldots,J.
\end{displaymath}
Then there exists $h_*>0$ such that for all $0<h \leq h_*$ and all 
$j=1,\ldots,J$ we have 
\begin{equation} \label{eq:c0qj}
\tfrac12 c_0 \leq q_j \leq 2C_0
\end{equation}
and
\begin{subequations} \label{eq:consist}
\begin{align}
\tfrac12 (q_j+q_{j+1}) & = | x_\rho(\rho_j) | + \mathcal O(h^2); \label{eq:consist1} \\
\tau_{j} + \tau_{j+1} & = 2 \, \tau(\rho_j) + \mathcal O(h^2); \label{eq:consist2} \\
\frac{\tau_{j+1} - \tau_j}{h} & = \tau_\rho(\rho_j) + \mathcal O(h^2). \label{eq:consist3}
\end{align}
\end{subequations}
\end{lemma}
\begin{proof} A Taylor expansion yields
\[
\delta x_{j+1} = \frac{x_{j+1} - x_j}{h}  =   x_\rho + \frac{h}{2} x_{\rho \rho} + \frac{h^2}{6} x_{\rho \rho \rho} + \mathcal O(h^3),
\]
where all the derivatives of $x$, and $\tau$, in this proof are evaluated at 
$\rho_j$. 
Hence
\begin{align*}
q_{j+1}^2 & = |x_\rho|^2 + h x_{\rho \rho} \cdot x_\rho 
+ \frac{h^2}{4} |x_{\rho \rho}|^2 + \frac{h^2}{3} x_{\rho \rho \rho} \cdot x_\rho
+ \mathcal O(h^3) \\ &
= |x_\rho|^2 \left( 1 + h \frac{x_{\rho \rho}}{| x_\rho |} \cdot \tau 
+ h^2 \left[ \tfrac14 \frac{ | x_{\rho \rho} |^2}{| x_\rho |^2} +
\tfrac13 \frac{x_{\rho \rho \rho}}{|x_\rho |} \cdot \tau \right]  
+ \mathcal O(h^3) \right),
\end{align*}
and with $\sqrt{1+\epsilon} = 1 + \frac12\epsilon - \frac18\epsilon^2 +
\mathcal{O}(\epsilon^3)$ it therefore follows that
\[
q_{j+1} = | x_\rho | \left( 1 + \frac{h}{2} \frac{x_{\rho \rho}}{| x_\rho |} \cdot \tau + h^2 \left[ \tfrac18 \frac{ | x_{\rho \rho} |^2}{| x_\rho |^2} +
\tfrac16 \frac{x_{\rho \rho \rho}}{|x_\rho |} \cdot \tau 
- \tfrac18 \frac{ (x_{\rho \rho} \cdot \tau)^2}{| x_\rho |^2}
\right] \right) 
+ \mathcal O(h^3).
\]
Moreover, since $\tau_{j+1} = \frac{1}{q_{j+1}} \delta x_{j+1}$ 
and $\frac1{1+\epsilon} = 1 - \epsilon + \epsilon^2 + \mathcal{O}(\epsilon^3)$,
we have that
\begin{align*}
\tau_{j+1} & = \tau + \frac{h}{2} \tau_\rho 
+ h^2 \left[ \tfrac16 \frac{x_{\rho \rho \rho}}{| x_\rho |} 
- \left[ \tfrac18 \frac{ | x_{\rho \rho} |^2}{| x_\rho |^2} +
\tfrac16 \frac{x_{\rho \rho \rho}}{|x_\rho |} \cdot \tau 
- \tfrac38 \frac{ (x_{\rho \rho} \cdot \tau)^2}{| x_\rho |^2} \right]\tau
-\tfrac14 (\frac{x_{\rho \rho}}{| x_\rho |} \cdot \tau) \frac{x_{\rho \rho}}{| x_\rho |}
 \right] + \mathcal O(h^3),
\end{align*}
where we used that $\displaystyle \frac{x_{\rho \rho}}{| x_\rho |} - (\frac{x_{\rho \rho}}{| x_\rho |} \cdot \tau) \tau = \bigl( \frac{x_\rho}{| x_\rho|} \bigr)_\rho = \tau_\rho$.
In a similar way one finds that
\begin{align*}
q_j & = | x_\rho | \left( 1 - \frac{h}{2} \frac{x_{\rho \rho}}{| x_\rho |} \cdot \tau 
+ h^2 \left[ \tfrac18 \frac{ | x_{\rho \rho} |^2}{| x_\rho |^2} 
+ \tfrac16 \frac{x_{\rho \rho \rho}}{|x_\rho |} \cdot \tau 
- \tfrac18 \frac{ (x_{\rho \rho} \cdot \tau)^2}{| x_\rho |^2} 
\right] \right) 
+ \mathcal O(h^3); \\
\tau_{j} & = \tau - \frac{h}{2} \tau_\rho 
+ h^2 \left[ \tfrac16 \frac{x_{\rho \rho \rho}}{| x_\rho |} 
- \left[ \tfrac18 \frac{ | x_{\rho \rho} |^2}{| x_\rho |^2} 
+ \tfrac16 \frac{x_{\rho \rho \rho}}{|x_\rho |} \cdot \tau 
- \tfrac38 \frac{ (x_{\rho \rho} \cdot \tau)^2}{| x_\rho |^2} 
\right]\tau 
-\tfrac14 (\frac{x_{\rho \rho}}{| x_\rho |} \cdot \tau) \frac{x_{\rho \rho}}{| x_\rho |} \right] + \mathcal O(h^3).
\end{align*}
{From} the above we infer that \eqref{eq:c0qj} holds provided that 
$0< h \leq h_*$. The estimates \eqref{eq:consist} also follow immediately. 
\end{proof}

In view of \eqref{eq:consist1}  
a natural semidiscrete finite difference approximation of \eqref{eq:hcsf}
 is now defined as follows.
Find $ x^h: \mathcal G^h \times [0,T] \to \bR^2$ such that 
\begin{subequations} \label{eq:fdsd}
\begin{alignat}{2} 
 \ddot{ x}^h_j +\beta \dot{x}^h_j & = \frac{2}{q^h_j + q^h_{j+1}} 
\frac{\tau^h_{j+1} - \tau^h_{j}}{h} - (\dot{ x}^h_j \cdot \dot{\theta}^h_j) \theta^h_j
\quad \text{in } [0,T], \quad &&
j = 1,\ldots,J; \label{eq:testtheta} \\
x^h_j(0) & = x_0(\rho_j), \
\dot{x}^h_j(0)= \mathcal{V}_0(\rho_j) \theta^{h,\perp}_j(0), 
&& j = 1,\ldots,J. \label{eq:discinit}
\end{alignat}
\end{subequations}
Standard ODE theory implies that the above system has a unique solution on some interval $[0,T_h)$. 
Let us begin by deriving discrete analogues of \eqref{eq:xnormal} and 
\eqref{eq:leequation}.

\begin{lemma} 
Let $x^h: \mathcal G^h \times [0,T_h) \to \bR^2$ be a solution of 
\eqref{eq:fdsd}. Then we have in $[0,T_h)$ and for all $j=1,\ldots,J$ that
\begin{subequations} 
\begin{align}
\dot{x}^h_j \cdot \theta^h_j & = 0; \label{eq:discle1} \\
\dot{q}^h_j 
+ \tfrac14 ( q^h_{j-1} + q^h_j) \bigl( \dot{x}^h_{j-1} \cdot \ddot{x}^h_{j-1} + \beta | \dot{x}^h_{j-1} |^2 \bigr) 
+ \tfrac14 ( q^h_j + q^h_{j+1}) \bigl( \dot{x}^h_j \cdot \ddot{x}^h_j + \beta | \dot{x}^h_j |^2 \bigr) 
&= 0. \label{eq:discle2}
\end{align}
\end{subequations}
\end{lemma}
\begin{proof}
It follows from \eqref{eq:testtheta}, \eqref{eq:fdthetaorthogonal} and the fact that $| \theta^h_j | =1$ that
\begin{equation*} 
(\dot{ x}^h_j \cdot \theta^h_j)_t = \ddot{x}^h_j \cdot \theta^h_j + \dot{x}^h_j \cdot \dot{\theta}^h_j = - \beta \, \dot{ x}^h_j \cdot \theta^h_j ,\quad j = 1,\ldots,J.
\end{equation*}
Since $\dot{ x}^h_j(0) \cdot \theta^h_j(0) = 0$ by \eqref{eq:discinit}, 
we deduce \eqref{eq:discle1}. In particular, $\dot{x}^h_j \cdot \tau^h_j = - \dot{x}^h_{j} \cdot \tau^h_{j+1}$
and hence
\begin{align*}
\dot{q}^h_j & = \frac{\dot{x}^h_j - \dot{x}^h_{j-1}}{h} \cdot \tau^h_j 
= - \tfrac12 \dot{x}^h_j \cdot \frac{ \tau^h_{j+1} - \tau^h_j}{h}
- \tfrac12 \dot{x}^h_{j-1} \cdot \frac{ \tau^h_j - \tau^h_{j-1}}{h} \\
& = -\tfrac14 ( q^h_j + q^h_{j+1}) \bigl( \dot{x}^h_j \cdot \ddot{x}^h_j + \beta | \dot{x}^h_j |^2 \bigr)
- \tfrac14 ( q^h_{j-1} + q^h_j) \bigl( \dot{x}^h_{j-1} \cdot \ddot{x}^h_{j-1} + \beta | \dot{x}^h_{j-1} |^2 \bigr),
\end{align*}
where the last equation is a consequence of \eqref{eq:testtheta} and 
\eqref{eq:discle1}. This proves \eqref{eq:discle2}.
\end{proof}

We also have the following discrete analogue of Remark~\ref{rem:normalflow}, 
where for simplicity we consider only the case $\beta = 0$.
\begin{lemma} \label{lem:stab}
Let $x^h: \mathcal G^h \times [0,T_h) \to \bR^2$ be a solution of 
\eqref{eq:fdsd} with $\beta = 0$. Then we have in $[0,T_h)$ that
\begin{equation} \label{eq:stab}
\tfrac12 \ddt h \sum_{j=1}^J \left[
\tfrac12 (q^h_j + q^h_{j+1}) | \dot x^h_j|^2 + 2 q^h_j \right] 
= \tfrac12 h \sum_{j=1}^J \dot q^h_j 
 \tfrac12 ( | \dot x^h_{j-1}|^2 + | \dot x^h_j|^2 )
,
\end{equation}
with $\dot q^h_j$ satisfying \eqref{eq:discle2}.
\end{lemma}
\begin{proof}
We compute, on noting \eqref{eq:testtheta} and \eqref{eq:discle1}, that
\begin{align}
\tfrac12 \ddt h \sum_{j=1}^J \tfrac12 (q^h_j + q^h_{j+1}) | \dot x^h_j|^2 
& =  \tfrac12 h \sum_{j=1}^J \tfrac12 (\dot q^h_j + \dot q^h_{j+1}) 
| \dot x^h_j|^2 
+ h \sum_{j=1}^J \tfrac12 (q^h_j + q^h_{j+1}) \dot x^h_j \cdot \ddot x^h_j  
\nonumber \\ & 
= \tfrac12 h \sum_{j=1}^J \dot q^h_j 
\tfrac12 ( | \dot x^h_{j-1}|^2 + | \dot x^h_j|^2 )
+ h \sum_{j=1}^J \dot x^h_j \cdot \frac{\tau^h_{j+1} - \tau^h_{j}}{h}
\nonumber \\ & 
= \tfrac12 h \sum_{j=1}^J \dot q^h_j 
\tfrac12 ( | \dot x^h_{j-1}|^2 + | \dot x^h_j|^2 )
- h \sum_{j=1}^J \frac{\dot x^h_j - \dot x^h_{j-1}}h \cdot \tau^h_{j}
\nonumber \\ & 
= \tfrac12 h \sum_{j=1}^J \dot q^h_j 
\tfrac12 ( | \dot x^h_{j-1}|^2 + | \dot x^h_j|^2 )
- \ddt h \sum_{j=1}^J q^h_j .
\label{eq:stabproof}
\end{align}
Combining \eqref{eq:stabproof} with \eqref{eq:discle2} then yields the desired
result \eqref{eq:stab}. 
\end{proof}

Observe that the right hand side of \eqref{eq:stab},
in view of \eqref{eq:discle2}, approximates the expression
\begin{align}
& - \tfrac12 \int_I |x_\rho| (x_t \cdot x_{tt} 
) |x_t|^2 \drho
=  - \tfrac12 \int_I |x_\rho| (\varkappa x_t \cdot \nu) |x_t|^2 \drho
,
\label{eq:stabc}
\end{align}
where we have noted \eqref{eq:abxttxss}, \eqref{eq:nu} and \eqref{eq:xnormal}.
As \eqref{eq:stabc} agrees with the right hand side in 
\eqref{eq:lemdte} with $\beta = 0$, recall again \eqref{eq:xnormal}, 
Lemma~\ref{lem:stab} can be viewed as a discrete analogue of 
Remark~\ref{rem:normalflow}. 

We stress that utilizing a suitable analogue of \eqref{eq:stab} 
will be at the heart of our error
analysis in Section~\ref{sec:ee}, below. In particular, 
$\dot x^h_j$ will be replaced
by the time derivative of the error between $x$ and $x^h$ at the point 
$\rho_j$.

Let us next consider the consistency errors for the scheme 
\eqref{eq:testtheta} and for the property \eqref{eq:discle2}.

\begin{lemma} \label{lem:consist} 
Let $x$ be the solution of \eqref{eq:hcsf}. Define
\begin{subequations} 
\begin{align}
R_j &:=  \ddot{ x}_j  + \beta \dot{x}_j- \frac{2}{q_j + q_{j+1}}   \frac{\tau_{j+1} - \tau_{j}}{h} +  (\dot{ x}_j \cdot \tau_t(\rho_j,\cdot)) \tau(\rho_j,\cdot); 
\label{eq:defRj} \\
\tilde R_j & := \dot{q}_j  + \tfrac14 ( q_{j-1}+ q_j) \bigl( \dot{x}_{j-1} \cdot \ddot{x}_{j-1} + \beta | \dot{x}_{j-1} |^2 \bigr)
 +\tfrac14 ( q_{j}+ q_{j+1}) \bigl( \dot{x}_{j} \cdot \ddot{x}_{j} + \beta | \dot{x}_j |^2 \bigr). \label{eq:qdot}
\end{align}
\end{subequations}
Then there exists a constant $C_1$ such that
\begin{equation} \label{eq:consist4}
\max_{j=1,\ldots,J} \bigl( |R_j(t) | + | \tilde R_j(t) | \bigr) \leq C_1 h^2,  \quad 
\ t \in [0,T].
\end{equation}
\end{lemma}
\begin{proof} 
The bound on $R_j$ is a direct consequence of Lemma~\ref{lem:prop1}. 
In order to analyze $\tilde R_j$ we deduce from \eqref{eq:consist3} that
$\tau_{j \pm 1} = \tau_j \pm h \tau_\rho(\rho_j,\cdot) + \mathcal O(h^2)$,
and hence by \eqref{eq:consist2}
\begin{displaymath}
\tau_j = \tfrac12 \frac{\tau_j+\tau_{j+1}}{2} 
+ \tfrac12 \frac{\tau_{j-1}+\tau_j}{2} + \mathcal O(h^2) 
= \tfrac12 \bigl( \tau(\rho_j,\cdot) + \tau(\rho_{j-1},\cdot)
\bigr) + \mathcal O(h^2).
\end{displaymath}
Combining this relation with the fact that 
\begin{displaymath}
\frac{\dot{x}_j - \dot{x}_{j-1}}{h} = \tfrac12 \bigl( x_{t \rho}(\rho_{j-1},\cdot)+ x_{t \rho}(\rho_j,\cdot) \bigr) + \mathcal O(h^2)
\end{displaymath}
we obtain
\begin{align*}
\dot{q}_j & = \frac{\dot{x}_j - \dot{x}_{j-1}}{h} \cdot \tau_j 
= \tfrac14 \bigl( x_{t \rho}(\rho_{j-1},\cdot)+ x_{t \rho}(\rho_j,\cdot) \bigr) \cdot
\bigl( \tau(\rho_{j-1},\cdot) + \tau(\rho_j,\cdot) \bigr) + \mathcal O(h^2) \\
& = \tfrac12 x_{t \rho}(\rho_{j-1},\cdot) \cdot \tau(\rho_{j-1},\cdot) 
+ \tfrac12 x_{t \rho}(\rho_j,\cdot) \cdot \tau(\rho_j,\cdot) \\ & \qquad
- \tfrac14 \bigl( x_{t \rho}(\rho_{j},\cdot)- x_{t \rho}(\rho_{j-1},\cdot) \bigr) \cdot \bigl( \tau(\rho_{j},\cdot) - \tau(\rho_{j-1},\cdot) \bigr) + \mathcal O(h^2) \\
& = \tfrac12 x_{t \rho}(\rho_{j-1},\cdot) \cdot \tau(\rho_{j-1},\cdot) 
+ \tfrac12 x_{t \rho}(\rho_j,\cdot) \cdot \tau(\rho_j,\cdot) + \mathcal O(h^2) \\
& = - \tfrac12 | x_\rho(\rho_{j-1},\cdot)|  \bigl( \dot x_{j-1} \cdot \ddot x_{j-1} + \beta | \dot{x}_{j-1} |^2 \bigr)
- \tfrac12 | x_\rho(\rho_{j},\cdot)| \bigl( \dot x_{j}\cdot \ddot x_{j} + \beta | \dot{x}_j |^2 \bigr)
+ \mathcal O(h^2),
\end{align*}
where we have used \eqref{eq:leequation}. Now the bound on $\tilde R_j$ follows with the help of \eqref{eq:consist1}.
\end{proof}

In view of Lemma \ref{lem:consist} we expect second order convergence for our scheme. As our  main result we prove
that this is indeed the case, where the error is measured in discrete integral norms that are natural for a
second order system of hyperbolic PDEs.
\begin{theorem} \label{thm:main} 
Suppose that Assumption~\ref{ass:x} is satisfied. 
Then there exists $h_0>0$ such that for $0 < h \leq h_0$
the problem \eqref{eq:fdsd}
has a unique solution $x^h: \mathcal G^h \times [0,T] \to \bR^2$ and the following error bounds hold:
\begin{equation} \label{eq:errbound}
\max_{0 \leq t \leq T} \bigl( \| x(t) - x^h(t) \|_{1,h} + \| \dot{x}(t) - \dot{x}^h(t) \|_{0,h} \bigr) \leq C h^2.
\end{equation}
\end{theorem}

Here, and throughout, $C$ denotes a generic positive constant independent of 
the mesh parameter $h$.

\setcounter{equation}{0}
\section{Proof of Theorem \ref{thm:main}} \label{sec:ee}

Let us abbreviate
\begin{displaymath}
x_j(t)=x(\rho_j,t), \; q_j(t) = | \delta x_j(t)|, \; \mbox{ and } \tau_{j}(t) = \frac{1}{q_{j}(t)} \delta x_{j}(t), \quad j=1,\ldots,J,
\end{displaymath}
where $x$ denotes the solution of \eqref{eq:hcsf}. Furthermore, we let
\begin{align} 
\hat T_h = \sup \Bigl\{ \hat t \in [0,T] & :
x^h \mbox{ solves } \eqref{eq:fdsd} \mbox{ on } [0,\hat t],
\text{ with } \tfrac14 c_0 \leq q^h_j(t) \leq 4 C_0 
\text{ and } \nonumber \\ & \
\max_{j=1,\ldots,J} \left(
| \tau_j(t) - \tau^h_j(t) | + | \dot{x}_j(t) - \dot{x}^h_j(t) | \right)
\leq h^{\frac54} \ \mbox{for } 0 \leq t \leq \hat t \Bigr\}. \label{eq:hatTh}
\end{align}
Here we have chosen the power $h^\frac54$ in the definition \eqref{eq:hatTh} 
as a convenient value between $1$ and $\frac32$, where the latter power of $h$
arises in the proof due to the application of an inverse inequality,
see \eqref{eq:h54} below.

Clearly, $\hat T_h>0$. In view of \eqref{eq:lereg} and Lemma~\ref{lem:prop1} 
we may assume that
\begin{subequations}
\begin{equation} \label{eq:tautau}
| \tau_j + \tau_{j+1} | \geq 1 
\end{equation}
and hence
\begin{equation} \label{eq:tauhtauh}
| \tau^h_j + \tau^h_{j+1} | \geq | \tau_j + \tau_{j+1} | - | \tau^h_j - \tau_j | - | \tau^h_{j+1} - \tau_{j+1} | \geq 1 - 2 h^{\frac{5}{4}} \geq \tfrac12,
\end{equation}
\end{subequations}
provided that $0<h \leq h_*$ is sufficiently small.
Thus $\theta^h_j(t)$ is well defined for $j=1,\ldots,J$ and 
$t \in [0,\hat T_h)$. Furthermore, we have:

\begin{lemma} \label{lem:discapriori1}
There exists $0<h_0 \leq h_*$ and a constant $C_2$, which only depends on $c_0, C_0$ and $\beta$,
such that for all $0<h \leq h_0$ and $0 \leq t < \hat T_h$
\begin{displaymath}
\max_{j=1,\ldots,J} \left(
| \dot{x}^h_j(t) | + | \dot{\theta}^h_j(t) | + | \ddot{x}^h_j(t) | + | \dot{q}^h_j(t) | \right)
\leq C_2 .
\end{displaymath}
\end{lemma}
\begin{proof} 
To begin, we deduce from \eqref{eq:lereg} and \eqref{eq:hatTh} that 
\[ | \dot{x}^h_j(t) | \leq | x_t(\rho_j,t) | + | (\dot{x}_j - \dot{x}^h_j)(t) | \leq C_0 + h^{\frac{5}{4}} \leq 2 C_0,
\]
provided that $0<h \leq h_0$ with $h_0$ sufficiently small.
Next, a straightforward calculation shows that $\displaystyle 
\dot{\tau}^h_j = \frac{1}{q^h_j} \bigl( \delta \dot{x}^h_j - ( \delta \dot{x}^h_j \cdot \tau^h_j) \tau^h_j \bigr)$ 
and hence, on noting \eqref{eq:hatTh}, \eqref{eq:delta-} and \eqref{eq:lereg},
it holds that 
\begin{align*}
| \dot{\tau}^h_j(t) | & \leq \frac{1}{q^h_j(t)} | \delta \dot{x}^h_j(t) | \leq \frac{4}{c_0} \bigl(| \delta (\dot{x}^h_j - \dot{x}_j)(t) | + | \delta \dot{x}_j(t) | \bigr) \\
& \leq \frac{8}{c_0} \frac{1}{h} \max_{1 \leq k \leq J} | (\dot{x}^h_k - \dot{x}_k)(t) | + \frac{4}{c_0} \max_{\rho \in I} | x_{t \rho}(\rho,t)| \leq \frac{8}{c_0} h^{\frac{1}{4}} + \frac{4}{c_0} C_0 \leq \frac{ 8 C_0}{c_0},
\end{align*}
provided that $0<h \leq h_0$ with $h_0$ sufficiently small. 
{From} this we deduce, on recalling \eqref{eq:fdthetah} and
\eqref{eq:tauhtauh}, that
\[
|\dot{\theta}^h_j(t)| \leq
\frac{|\dot\tau^h_j(t) + \dot\tau^h_{j+1}(t)|}{|\tau^h_j(t) + \tau^h_{j+1}(t)|}
\leq 2 \frac{ 16 C_0}{c_0} = \frac{ 32 C_0}{c_0}.
\]
In order to bound $\ddot{x}^h_j$, we first use \eqref{eq:consist3}, 
\eqref{eq:lereg} and \eqref{eq:hatTh} to show that
\begin{displaymath}
\left| \frac{\tau^h_{j+1}(t) - \tau^h_j(t)}{h} \right| 
\leq \left| \frac{\tau_{j+1}(t) - \tau_j(t)}{h} \right| 
+ \frac{2}{h} \max_{k=1,\ldots,J} | \tau_k(t) - \tau^h_k(t) | 
\leq 2 C_0 + 2 h^{\frac{1}{4}} \leq 3 C_0,
\end{displaymath}
provided that $0<h \leq h_0$ with $h_0\leq h_*$ sufficiently small. 
If we combine this estimate with \eqref{eq:testtheta}, \eqref{eq:hatTh} and
the previously derived bounds on $\dot{x}^h_j$ and $\dot{\theta}^h_j$,
we obtain
\begin{align*}
|\ddot{x}^h_j(t)| & \leq \beta |\dot{x}^h_j(t)| + \frac2{q^h_j(t)+q^h_{j+1}(t)}
\left|\frac{\tau^h_{j+1}(t) - \tau^h_{j}(t)}{h}\right| + 
|\dot{ x}^h_j(t)| \, |\dot{\theta}^h_j(t)| \\ &
\leq 2 \beta C_0 + \frac{4}{c_0} 3 C_0 + 2 C_0 \frac{ 32 C_0 }{c_0} 
= 2 \beta C_0 + \frac{12 C_0}{c_0} + \frac{ 64 C_0^2 }{c_0}.
\end{align*}
Finally, the bound on $\dot{q}^h_j$ is a
consequence of \eqref{eq:discle2} and \eqref{eq:hatTh} using now in addition the bound on $\ddot{x}^h_j$.  
\end{proof}

Let us introduce the error $e_j(t):=x_j(t)-x^h_j(t)$. 
We infer from \eqref{eq:defRj} and \eqref{eq:testtheta} that
\begin{align}
&  \ddot{e}_j + \beta \dot{e}_j - \frac{2}{q^h_j+ q^h_{j+1}} \frac{(\tau_{j+1} - \tau^h_{j+1})- (\tau_j - \tau^h_j)}{h} 
\nonumber \\ & \quad 
=  \bigl( \dot{ x}^h_j \cdot (\dot{\theta}^h_j - \tau_t(\rho_j,\cdot)  \bigr) \tau(\rho_j,\cdot) + 
( \dot{x}^h_j \cdot \dot{\theta}^h_j) \bigl( \theta^h_j - \tau(\rho_j,\cdot) \bigr) -  (\dot{ e}_j \cdot \tau_t(\rho_j,\cdot)) \tau(\rho_j,\cdot)
\nonumber \\ & \qquad
+ 2 \frac{  (q^h_j - q_j) + (q^h_{j+1}-  q_{j+1})}{(q_j+q_{j+1}) (q^h_j+q^h_{j+1})} \frac{\tau_{j+1} - \tau_j}{h} + R_j
\nonumber  \\ & \quad
=: \sum_{k=1}^5 T^k_j. \label{eq:err1}
\end{align}
Taking the scalar product with $\frac{h}{2}(q^h_j+ q^h_{j+1}) \dot{e}_j$, summing over $j=1,\ldots,J$ and recalling 
Lemma~\ref{lem:discapriori1} yields
\begin{align}
& \tfrac12 h \ddt \sum_{j=1}^J \tfrac12 (q^h_j+ q^h_{j+1}) | \dot{e}_j |^2  + \beta h \sum_{j=1}^J  \tfrac12 (q^h_j+ q^h_{j+1}) | \dot{e}_j |^2 
- h \sum_{j=1}^J \frac{(\tau_{j+1} - \tau^h_{j+1})- (\tau_j - \tau^h_j)}{h} 
\cdot \dot{e}_j 
\nonumber \\ & \quad
= \tfrac12 h \sum_{j=1}^J \tfrac12 (\dot{q}^h_j+ \dot{q}^h_{j+1}) |\dot e_j|^2 
+ \sum_{k=1}^5 h \sum_{j=1}^J  \tfrac12 (q^h_j+ q^h_{j+1})  T^k_j \cdot \dot{e}_j \nonumber \\ &  \quad
\leq C h \sum_{j=1}^J | \dot{e}_j |^2 
+ h \sum_{k=1}^5 \sum_{j=1}^J  \tfrac12 (q^h_j+ q^h_{j+1}) T^k_j \cdot \dot{e}_j. \label{eq:err2}
\end{align}
While the above relation already provides us with some control on $\dot{e}_j$, the treatment of the elliptic part
is more difficult. This is a consequence of the fact that the operator $\frac{1}{| x_\rho |} \bigl( \frac{x_\rho}{| x_\rho |} \bigr)_\rho$ is degenerate in tangential
direction. It is therefore not possible to directly control $\delta e_j$, which we split instead as follows:
\begin{equation} \label{eq:firstderiv}
\delta e_j = \delta x_j - \delta x^h_j = q_j (\tau_j - \tau^h_j) + (q_j - q^h_j) \tau^h_j.
\end{equation}
In the next step we will gain control on the difference of the tangents from the third term on the left hand side of \eqref{eq:err2}.
To do so, we essentially adapt arguments from \cite[Section 5]{Dziuk99} developed for a finite element approach to the
curve shortening flow. To begin, using summation by parts together with the fact that 
$\delta x^h_j = q^h_j \tau^h_j$, we derive
\begin{align*}
& - h \sum_{j=1}^J \frac{(\tau_{j+1} - \tau^h_{j+1})- (\tau_j - \tau^h_j)}{h} 
\cdot \dot{e}_j \\ & \quad
= h \sum_{j=1}^J (\tau_j - \tau^h_j) \cdot \frac{\dot{e}_j - \dot{e}_{j-1}}{h} 
= h \sum_{j=1}^J (\tau_j - \tau^h_j) \cdot
\bigl( \delta \dot{x}_j - \delta \dot{x}^h_j \bigr) \\ & \quad
= h \sum_{j=1}^J \bigl( \tau^h_j \cdot \delta \dot{x}^h_j - \tau_j \cdot
\delta \dot{x}^h_j \bigr) 
+ h \sum_{j=1}^J (\tau_j - \tau^h_j)  \cdot  \delta \dot{x}_j  \\ & \quad
= h \ddt \sum_{j=1}^J \bigl( q^h_j - \tau_j \cdot \delta x^h_j \bigr)
+ h \sum_{j=1}^J \dot{\tau}_j \cdot  \delta x^h_j 
+ h \sum_{j=1}^J (\tau_j - \tau^h_j)  \cdot  \delta \dot{x}_j \\ & \quad
= h \ddt \sum_{j=1}^J q^h_j ( 1 - \tau_j \cdot \tau^h_j) 
+ h \sum_{j=1}^J \frac{q^h_j}{q_j}  \bigl(  \delta \dot{x}_j  -
( \delta \dot{x}_j \cdot \tau_j) \tau_j \bigr) \cdot \tau^h_j 
+ h \sum_{j=1}^J (\tau_j - \tau^h_j)  \cdot \delta \dot{x}_j  \\ & \quad
= \tfrac12 h \ddt \sum_{j=1}^J q^h_j | \tau_j - \tau^h_j |^2 
+ h \sum_{j=1}^J \left( 
\frac{q_j- q^h_j}{q_j} \; \delta \dot{x}_j \cdot ( \tau_j - \tau^h_j)  
 + \tfrac12 \frac{q^h_j}{q_j} \left( \delta \dot{x}_j \cdot \tau_j \right)
 | \tau_j - \tau^h_j |^2 \right).
\end{align*}
If we insert the above relation into \eqref{eq:err2}, note that $\beta \geq 0$ and apply a
Cauchy--Schwarz inequality together with Lemma~\ref{lem:prop1} and 
\eqref{eq:hatTh}, we obtain
\begin{align}
&\tfrac12 h\ddt \sum_{j=1}^J \bigl( \tfrac12 (q^h_j+ q^h_{j+1}) | \dot{e}_j |^2 
+ q^h_j | \tau_j - \tau^h_j |^2 \bigr) \nonumber \\ & \qquad
\leq C h \sum_{j=1}^J \bigl( | \dot{e}_j |^2 + (q_j - q^h_j)^2 
+ |\tau_j - \tau^h_j |^2 \bigr) 
+ h \sum_{k=1}^5 \sum_{j=1}^J  \tfrac12 (q^h_j+ q^h_{j+1})  T^k_j \cdot \dot{e}_j. \label{eq:err3}
\end{align}
Let us next consider the terms involving $T^k_j,k=1,\ldots,5$. To begin, note that \eqref{eq:xnormal} and \eqref{eq:discle1} imply
\begin{equation} \label{eq:taudote}
\tau(\rho_j,\cdot) \cdot \dot{e}_j = \tau(\rho_j,\cdot) \cdot (\dot{x}_j - \dot{x}^h_j) = - \tau(\rho_j,\cdot) \cdot \dot{x}^h_j = \dot{x}^h_j \cdot ( \theta^h_j - \tau(\rho_j,\cdot)).
\end{equation}
Therefore the definition of $T^1_j$ in \eqref{eq:err1} yields that
\begin{align*}
& h \sum_{j=1}^J  \tfrac12 (q^h_j+ q^h_{j+1})  T^1_j \cdot \dot{e}_j = 
h \sum_{j=1}^J \tfrac12 (q^h_j + q^h_{j+1}) 
\bigl( \dot{x}^h_j \cdot (\dot{\theta}^h_j - \tau_t(\rho_j,\cdot)) \bigr) 
\bigl( \dot{x}^h_j \cdot ( \theta^h_j - \tau(\rho_j,\cdot)) \bigr) \\ & \quad
= \tfrac12 h \ddt \sum_{j=1}^J \tfrac12 (q^h_j + q^h_{j+1}) 
\bigl( \dot{ x}^h_j \cdot (\theta^h_j - \tau(\rho_j,\cdot))  \bigr)^2 - 
\tfrac12 h \sum_{j=1}^J \tfrac12 (\dot{q}^h_j + \dot{q}^h_{j+1}) 
\bigl( \dot{x}^h_j \cdot (\theta^h_j - \tau(\rho_j,\cdot))  \bigr)^2 
\\ & \qquad
- h \sum_{j=1}^J \tfrac12 (q^h_j + q^h_{j+1}) 
\bigl( \dot{ x}^h_j \cdot (\theta^h_j - \tau(\rho_j,\cdot))  \bigr) 
\bigl( \ddot{x}^h_j \cdot ( \theta^h_j - \tau(\rho_j,\cdot)) \bigr) \\ & \quad
\leq \tfrac12 h \ddt \sum_{j=1}^J \tfrac12 (q^h_j + q^h_{j+1}) 
\bigl( \tau(\rho_j,\cdot) \cdot \dot{e}_j  \bigr)^2
+ C h \sum_{j=1}^J | \theta^h_j - \tau(\rho_j,\cdot) |^2,
\end{align*}
where in the last step we have used \eqref{eq:taudote}, 
as well as Lemma~\ref{lem:discapriori1}. 
Moreover, we have from \eqref{eq:consist2} that
\begin{align*}
| \theta^h_j - \tau(\rho_j,\cdot) | & \leq 
\left| \frac{\tau^h_j+\tau^h_{j+1} }{| \tau^h_j+\tau^h_{j+1} |} 
- \frac{\tau_j+\tau_{j+1} }{| \tau_j+\tau_{j+1} |}  \right| 
+ \left| \frac{\tau_j+\tau_{j+1} }{| \tau_j+\tau_{j+1} |} 
-  \tau(\rho_j,\cdot) \right| \nonumber \\ & 
\leq \frac{2}{| \tau_j + \tau_{j+1} |} 
\bigl( | \tau_j - \tau^h_j| + | \tau_{j+1} - \tau^h_{j+1} | \bigr) 
+ Ch^2, 
\end{align*}
so that with \eqref{eq:tautau} 
\begin{equation} \label{eq:thetaminustau2}
|\theta^h_j - \tau(\rho_j,\cdot) |^2 \leq C (| \tau_j - \tau^h_j |^2 
+ | \tau_{j+1} - \tau^h_{j+1} |^2 ) + Ch^4.
\end{equation}
In particular, it follows that
\begin{equation} \label{eq:t1}
h \sum_{j=1}^J  \tfrac12 (q^h_j+ q^h_{j+1})  T^1_j \cdot \dot{e}_j 
\leq \tfrac12 h \ddt \sum_{j=1}^J \tfrac12 (q^h_j + q^h_{j+1}) 
\bigl( \tau(\rho_j,\cdot) \cdot \dot{e}_j  \bigr)^2
+ C h \sum_{j=1}^J | \tau_j - \tau^h_j |^2 + Ch^4.
\end{equation}
Next, we deduce with the help of Lemma~\ref{lem:discapriori1}, \eqref{eq:hatTh}
and \eqref{eq:thetaminustau2} that
\begin{equation} \label{eq:t2}
 h \sum_{j=1}^J \tfrac12 (q^h_j+ q^h_{j+1})  T^2_j \cdot \dot{e}_j \leq C h \sum_{j=1}^J | \theta^h_j - \tau(\rho_j,\cdot) | \, | \dot{e}_j | \leq
 C h \sum_{j=1}^J \bigl( | \tau_j-\tau^h_j |^2 + | \dot{e}_j |^2 \bigr) + C h^4,
\end{equation}
while in view of \eqref{eq:consist4}
\begin{equation} \label{eq:t3t5}
 h \sum_{j=1}^J \tfrac12 (q^h_j+ q^h_{j+1}) (T^3_j+ T^5_j) \cdot \dot{e}_j \leq C h \sum_{j=1}^J | \dot{e}_j |^2 + C h^4.
 \end{equation}
Finally, with the help of \eqref{eq:consist3} and \eqref{eq:lereg} we can bound
\begin{align} \label{eq:t4}
h \sum_{j=1}^J \tfrac12 (q^h_j+ q^h_{j+1}) T^4_j \cdot \dot{e}_j & 
\leq C h \sum_{j=1}^J \bigl( | q_j - q^h_j| + | q_{j+1}-q^h_{j+1} | \bigr) 
\frac{|\tau_{j+1} - \tau_j|}{h}
| \dot{e}_j | \nonumber \\ &
\leq  C h \sum_{j=1}^J \bigl( (q_j - q^h_j)^2 + | \dot{e}_j |^2 \bigr)
+C h^4.
\end{align}

If we insert \eqref{eq:t1}, \eqref{eq:t2}, \eqref{eq:t3t5} and \eqref{eq:t4} 
into the estimate \eqref{eq:err3} we obtain, 
upon subtracting the first term on the right hand side of
\eqref{eq:t1} from both sides of the inequality
and on noting $| \dot{e}_j |^2 - (\dot{e}_j \cdot \tau)^2 = 
(\dot{e}_j \cdot \nu)^2$, that
\begin{align}
& \tfrac12 h \ddt \sum_{j=1}^J \left( \tfrac12 (q^h_j+ q^h_{j+1}) 
\bigl( \dot{e}_j \cdot \nu(\rho_j,\cdot) \bigr)^2 
+ q^h_j | \tau_j - \tau^h_j |^2 \right) \nonumber \\ & \quad
\leq C h \sum_{j=1}^J 
\left( | \dot{e}_j |^2 + (q_j - q^h_j)^2 + |\tau_j - \tau^h_j |^2 \right) 
+ C h^4. \label{eq:err4}
\end{align}
Using \eqref{eq:taudote}, \eqref{eq:thetaminustau2} 
and Lemma~\ref{lem:discapriori1}, we have
\begin{align} \label{eq:normet}
h \sum_{j=1}^J | \dot{e}_j |^2 &
= h \sum_{j=1}^J \bigl( ( \dot{e}_j \cdot \tau(\rho_j,\cdot))^2 
+ (\dot{e}_j \cdot \nu(\rho_j,\cdot))^2 \bigr) \nonumber \\ & 
= h \sum_{j=1}^J \left(\dot{x}^h_j \cdot ( \theta^h_j - \tau(\rho_j,\cdot))
\right)^2
+ h \sum_{j=1}^J (\dot{e}_j \cdot(\nu(\rho_j,\cdot))^2 \nonumber \\ &
\leq C h \sum_{j=1}^J | \tau_j - \tau^h_j |^2 + C h^4 
+ h \sum_{j=1}^J (\dot{e}_j \cdot(\nu(\rho_j,\cdot))^2.
\end{align}
If we insert \eqref{eq:normet} into \eqref{eq:err4}  we find 
\begin{equation} \label{eq:err5}
\phi_h'(t) \leq C_3 \bigl( h^4 + \phi_h(t) + \psi_h(t) \bigr),
\end{equation}
where we have abbreviated 
\begin{equation} \label{eq:phipsi}
\phi_h(t):= h \sum_{j=1}^J \bigl( \tfrac12 (q^h_j+ q^h_{j+1}) \bigl( \dot{e}_j \cdot \nu(\rho_j,\cdot) \bigr)^2 + q^h_j | \tau_j - \tau^h_j |^2 \bigr), \quad
\psi_h(t):= h \sum_{j=1}^J (q_j - q^h_j)^2 ,
\end{equation}
and noted \eqref{eq:hatTh}.

It remains to bound the function $\psi_h$, which controls the second part in \eqref{eq:firstderiv}.
To do so we combine \eqref{eq:qdot} and \eqref{eq:discle2} and obtain
\begin{align} \label{eq:qjqhj}
\dot{q}_j - \dot{q}^h_j 
& = - \tfrac14 ( q_{j-1}+ q_{j}) 
\bigl( \dot{x}_{j-1} \cdot \ddot{x}_{j-1}- \dot{x}^h_{j-1} \cdot 
\ddot{x}^h_{j-1} \bigr)
- \tfrac14 ( q_j+ q_{j+1}) 
\bigl( \dot{x}_j \cdot \ddot{x}_j - \dot{x}^h_j \cdot \ddot{x}^h_j \bigr) 
\nonumber \\ & \quad
+ \tfrac14 \bigl( (q^h_{j-1}-q_{j-1}) + (q^h_{j} -q_{j}) \bigr) 
\dot{x}^h_{j-1} \cdot \ddot{x}^h_{j-1} 
+ \tfrac14 \bigl( (q^h_j-q_j) + (q^h_{j+1} -q_{j+1}) \bigr) 
\dot{x}^h_j \cdot \ddot{x}^h_j \nonumber \\ &  \quad
 - \tfrac14\beta ( q_{j-1}+ q_{j}) \bigl( | \dot{x}_{j-1} |^2 - | \dot{x}^h_{j-1} |^2 \bigr) - 
    \tfrac14\beta ( q_{j}+ q_{j+1}) \bigl( | \dot{x}_{j} |^2 - | \dot{x}^h_{j} |^2 \bigr)   \nonumber \\ & \quad 
    + \tfrac14\beta \bigl( (q^h_{j-1}-q_{j-1}) + (q^h_{j} -q_{j}) \bigr) | \dot{x}^h_{j-1} |^2 
    +  \tfrac14\beta \bigl( (q^h_{j}-q_{j}) + (q^h_{j+1} -q_{j+1}) \bigr) | \dot{x}^h_{j} |^2+ \tilde R_j \nonumber \\ &
= - \tfrac18 \partial_t \Bigl( ( q_{j-1}+ q_{j}) 
\bigl( | \dot{x}_{j-1} |^2 - | \dot{x}^h_{j-1} |^2  \bigr) \Bigr) 
- \tfrac18 \partial_t \Bigl( ( q_{j}+ q_{j+1}) 
\bigl( | \dot{x}_j |^2 - | \dot{x}^h_j |^2 \Bigr) \nonumber \\ & \quad
+ \tfrac18 \bigl( \dot{q}_{j-1}+ \dot{q}_j \bigr) 
\bigl( | \dot{x}_{j-1} |^2 - | \dot{x}^h_{j-1} |^2 \bigr) 
+ \tfrac18 \bigl( \dot{q}_{j}+ \dot{q}_{j+1} \bigr) 
\bigl( | \dot{x}_{j} |^2 - | \dot{x}^h_{j} |^2 \bigr) \nonumber \\ & \quad
+ \tfrac14 \bigl( (q^h_{j-1}-q_{j-1}) + (q^h_{j} -q_{j}) \bigr) 
\dot{x}^h_{j-1} \cdot \ddot{x}^h_{j-1} 
+ \tfrac14 \bigl( (q^h_j-q_j) + (q^h_{j+1} -q_{j+1}) \bigr) 
\dot{x}^h_j \cdot \ddot{x}^h_j  \nonumber \\ & \quad
 - \tfrac14\beta ( q_{j-1}+ q_{j}) \bigl( | \dot{x}_{j-1} |^2 - | \dot{x}^h_{j-1} |^2 \bigr) - 
    \tfrac14\beta ( q_{j}+ q_{j+1}) \bigl( | \dot{x}_{j} |^2 - | \dot{x}^h_{j} |^2 \bigr)  \nonumber \\ & \quad
    + \tfrac14\beta \bigl( (q^h_{j-1}-q_{j-1}) + (q^h_{j} -q_{j}) \bigr) | \dot{x}^h_{j-1} |^2 
    +  \tfrac14\beta \bigl( (q^h_{j}-q_{j}) + (q^h_{j+1} -q_{j+1}) \bigr) | \dot{x}^h_{j} |^2+ \tilde R_j .
\end{align}
Recalling \eqref{eq:init} and \eqref{eq:discinit} we infer that
$q^h_j(0) = q_j(0)$ as well as
\begin{equation} \label{eq:dotej0}
|\dot{e}_j(0)| = 
| \dot{x}_j(0) - \dot{x}^h_j(0) | = \left| \mathcal V_0(\rho_j) \Bigl( \tau(\rho_j,0) - \frac{\tau_j(0) + \tau_{j+1}(0)}{| \tau_j(0) + \tau_{j+1}(0) |} \Bigr)^\perp 
\right| \leq C h^2
\end{equation}
where we also made use of \eqref{eq:consist2}.
Thus we obtain after integrating \eqref{eq:qjqhj} in time, on noting
that $||a|^2 - |b|^2| \leq (|a|+|b|)|a-b|$ and on taking into account 
Lemma~\ref{lem:discapriori1} and \eqref{eq:consist4}, that
\begin{align*}
| q_j (t)- q^h_j(t) | & 
\leq C \bigl( | \dot{e}_{j-1}(t) | + | \dot{e}_j(t) | \bigr) 
+ C \int_0^t | \dot{e}_{j-1}(u) | + | \dot{e}_j(u) | \du 
\\ & \quad 
+ C \int_0^t | (q_{j-1}-q^h_{j-1})(u) | + | (q_j-q^h_j)(u) | 
+ | (q_{j+1}-q^h_{j+1})(u) | \du + C h^2 \\ & 
\leq C \bigl( | \dot{e}_{j-1}(t) | + | \dot{e}_j(t) | \bigr) 
+ C \Bigl(  \int_0^t | \dot{e}_{j-1}(u) |^2+ | \dot{e}_j(u) |^2 \du 
\Bigr)^{\frac{1}{2}} \\ & \quad 
+ C \Bigl( \int_0^t | (q_{j-1}-q^h_{j-1})(u) |^2 + | (q_j-q^h_j)(u) |^2 
+ | (q_{j+1}-q^h_{j+1})(u) |^2 \du \Bigr)^{\frac{1}{2}} + C h^2.
\end{align*}
Taking the square and summing over $j$ yields 
\begin{displaymath}
h \sum_{j=1}^J (q_j-q^h_j)^2(t) 
\leq C h \sum_{j=1}^J | \dot{e}_j(t) |^2 
+ C \int_0^t h \sum_{j=1}^J | \dot{e}_j(u) |^2  \du 
+ C h \int_0^t (q_j-q^h_j)^2(u) \du + Ch^4,
\end{displaymath}
which together with \eqref{eq:normet} and \eqref{eq:hatTh} implies
\begin{equation} \label{eq:err6}
\psi_h(t)  \leq C_4 \Bigl( \phi_h(t) + \int_0^t \bigl( \phi_h(u) + \psi_h(u) \bigr) \du + h^4 \Bigr).
\end{equation}
If we multiply \eqref{eq:err5} by $2 C_4$, integrate with respect to time and 
combine the result with \eqref{eq:err6}, we obtain, on noting
from \eqref{eq:dotej0} that $\phi_h(0) \leq C_5 h^4$, that
\begin{displaymath}
C_4 \phi_h(t) + \psi_h(t) \leq (2 C_4 (C_3 + C_5) + C_4) 
\bigl( h^4 + \int_0^t \phi_h(u) + \psi_h(u) \du \bigr),
\end{displaymath}
from which we deduce with the help of Gronwall's lemma that
\begin{equation} \label{eq:err7}
\phi_h(t) + \psi_h(t) \leq C h^4, \quad 0 \leq t < \hat T_h.
\end{equation}
In particular, we have for $j=1,\ldots,J$ and $0  \leq t < \hat T_h$ that
\begin{equation} \label{eq:h54}
| (\tau_j - \tau^h_j)(t) | 
\leq h^{-\frac12} \Bigl( h \sum_{k=1}^J | (\tau_k - \tau^h_k)(t)| 
\Bigr)^{\frac12} 
\leq C h^{-\frac{1}{2}} \sqrt{\phi_h(t)} 
\leq C h^{\frac{3}{2}} \leq \tfrac12 h^{\frac{5}{4}},
\end{equation}
provided that $0 < h \leq h_0$ and $h_0$ is chosen smaller if necessary. 
In a similar way, on combining \eqref{eq:err7}, \eqref{eq:phipsi},
\eqref{eq:normet}, \eqref{eq:c0qj} and \eqref{eq:lereg}, we obtain that  
\begin{displaymath}
| (\dot{x}_j - \dot{x}^h_j)(t) | \leq \tfrac12 h^{\frac54}, 
\quad \tfrac13 c_0 \leq q^h_j(t) \leq 3 C_0, 
\quad j=1,\ldots,J, \ 0 \leq t < \hat T_h.
\end{displaymath}
If $\hat T_h < T$ one could therefore continue the discrete solution to an 
interval $[0, \hat T_h + \epsilon]$, for some $\epsilon>0$, such that
$\tfrac14 c_0 \leq q^h_j(t) \leq 4 C_0$, 
$|\tau_j(t) - \tau^h_j(t)| + |\dot{x}_j(t) - \dot{x}^h_j(t) |\leq h^{\frac54}$
for all $j=1,\ldots,J$ and $0 \leq t \leq \hat T_h + \epsilon$,
contradicting the definition of $\hat T_h$. 
Thus, $\hat T_h = T$. Finally, the bounds
\eqref{eq:errbound} follow from \eqref{eq:err7}, the definitions of $\phi_h$ and $\psi_h$, \eqref{eq:normet} and \eqref{eq:firstderiv}.

\setcounter{equation}{0}
\section{Numerical results} \label{sec:nr}

\subsection{Fully discrete scheme}

For the numerical simulations presented in this section,
we consider the following fully discrete approximation of \eqref{eq:fdsd}, 
where in order to discretize in time, we let $t_m=m\Delta t$, $m=0,\ldots,M$, 
with the uniform time step $\Delta t = \frac TM >0$. We will approximate 
$x^h(t_m)$ by the grid function $x^m : \mathcal{G}^h \to \bR^2$.
Analogously to \eqref{eq:fdthetah}, we define $\theta^m_j$, $j=1,\ldots,J$, 
in terms of $x^m$, and similarly for $q^m_j$ and $\tau^m_j$.
Then, given suitable initial data $x^0, x^{-1}: \mathcal G^h \to \bR^2$,
for $m=0,\ldots,M-1$ we find $x^{m+1}: \mathcal G^h \to \bR^2$ such 
that, for $j = 1,\ldots,J$, 
\begin{align} \label{eq:fdfdimpl}
& \tfrac12 (q^m_j + q^m_{j+1}) h
\frac{ x^{m+1}_j - 2  x^m_j +  x^{m-1}_j}{(\Delta t)^2} 
+ \tfrac14 \beta (q^m_j + q^m_{j+1}) h
\frac{ x^{m+1}_j - x^{m-1}_j}{\Delta t} 
\nonumber \\ & \quad 
= 
\frac{1}{2q^m_{j+1}} \delta x^{m+1}_{j+1} - \frac{1}{2q^m_j} \delta x^{m+1}_j
+\frac{1}{2q^m_{j+1}} \delta x^{m-1}_{j+1} - \frac{1}{2q^m_j} \delta x^{m-1}_j
\nonumber \\ & \qquad 
- \tfrac12 (q^m_j + q^m_{j+1}) h
(\frac{x^m_j -  x^{m-1}_j}{\Delta t} \cdot \frac{\theta^m_j - \theta^{m-1}_j}
{\Delta t}) \theta^m_j. 
\end{align}
Observe that we have chosen a linear discretization, that is
analogous to a
mass-lumped finite element approximation of \eqref{eq:abxttxss}, 
which uses
a semi-implicit approximation of $\frac{1}{| x_\rho |}
\bigl(\frac{x_\rho}{| x_\rho |} \bigr)_\rho$ in the spirit of the discretizations
proposed for the linear wave equation in e.g.\ \cite[\S2.7]{GrossmannR07}.
We remark that in contrast to the semidiscrete setting, recall 
\eqref{eq:discle1}, it does not appear possible to prove a fully discrete 
analogue of the crucial normal flow property \eqref{eq:xnormal} 
for the fully discrete scheme \eqref{eq:fdfdimpl}.
 
In order to derive suitable initial data for \eqref{eq:fdfdimpl}, 
we observe that the solution to \eqref{eq:hcsf} satisfies the Taylor expansion
\begin{align} \label{eq:x0Taylor}
x(\cdot, \Delta t) & = x + \Delta t x_t 
+ \tfrac12 (\Delta t)^2 x_{tt} + \mathcal{O}((\Delta t)^3)
\nonumber \\ &
= x + \Delta t \mathcal{V}_0 \nu 
+ \tfrac12 (\Delta t)^2 \left[ \frac{1}{| x_\rho |}
\bigl(\frac{x_\rho}{| x_\rho |} \bigr)_\rho 
- ( \mathcal{V}_0 \nu \cdot \tau_t) \tau - \beta  \mathcal{V}_0 \nu \right] 
+ \mathcal{O}((\Delta t)^3)
\nonumber \\ &
= x + \Delta t \mathcal{V}_0 \nu 
+ \tfrac12 (\Delta t)^2 \left[ \frac{1}{| x_\rho |}
\bigl(\frac{x_\rho}{| x_\rho |} \bigr)_\rho 
- \frac{1}{|  x_\rho |} \mathcal V_0 \mathcal V_{0,\rho} \tau - \beta  \mathcal{V}_0 \nu \right] 
+ \mathcal{O}((\Delta t)^3),
\end{align}
where on the right hand side we always evaluate $x$, $\tau$, $\nu$,
$\mathcal{V}_0$ and their derivatives at $(\cdot,0)$. 
Note in particular that in the last step we used that
\begin{displaymath}
\tau_t \cdot \nu = \frac{x_{t, \rho}}{| x_\rho |} \cdot \nu = \frac{1}{| x_\rho|} (\mathcal V_0 \nu)_\rho \cdot \nu =
\frac{1}{| x_\rho |} \mathcal V_{0,\rho}.
\end{displaymath}
Inspired by \eqref{eq:x0Taylor} we choose as initial data
\begin{align*} 
x^0_j & = x_0(\rho_j)\ \text{ and } \ \\
x^{-1}_j & = x^0_j - \Delta t \mathcal{V}_0(\rho_j) \theta^{0,\perp}_j \nonumber \\ & \qquad\
+ \tfrac12 (\Delta t)^2 \left[ 
\frac{2}{q^0_j + q^0_{j+1}} \left( \frac1h \left(
\frac{\delta x^{0}_{j+1}}{q^0_{j+1}} - \frac{\delta x^{0}_j}{q^0_j} 
\right) 
- \mathcal V_0(\rho_j) \mathcal V_{0,\rho}(\rho_j) \theta^0_j \right)
 - \beta V_0 \theta^{0,\perp}_j \right],
\end{align*}
for $j=1,\ldots,J$.

We stress that all our presented numerical experiments fall within the
scope of our main result, Theorem~\ref{thm:main}. However, for nonconvex
initial data, and for convex initial data with an initial normal velocity 
$\mathcal{V}_0$ such that $\max_I \mathcal{V}_0 > 0$, a rigorous
existence and regularity theory for the underlying PDE appears to be 
still lacking.

\subsection{Convergence experiment}

Our first set of numerical experiments is for the evolution of an initially
circular curve when $\beta=0$. 
It can be shown that a family of circles with 
radius $r(t)$ is a solution to \eqref{eq:Gurtinbeta0} with 
$\mathcal V_\Gamma\!\mid_{t=0} = V_0 \in \bR$
for the initial outer normal velocity, if
\begin{equation*} 
\ddot{r}(t) = - \frac1{r(t)} \quad \text{in } (0,T],\quad
r(0) = r_0,\ \dot{r}(0) = V_0.
\end{equation*}  
Upon integration we obtain that
\[
\tfrac12 (\dot{r}(t))^2 = \ln r_0 - \ln r(t) + \tfrac12 V_0^2
= \ln \frac{r_0}{r(t)} + \tfrac12 V_0^2.
\]
Hence
\begin{equation*} 
\dot{r}(t) = \pm \sqrt{2 \ln \frac{r_0}{r(t)} + V_0^2 },
\end{equation*}
which means that if $V_0 > 0$, then $r(t)$ will at first increase until it hits
a maximum, where $2\ln \frac{r_0}{r(t)} + V_0^2 = 0$, after which it
will decrease and shrink to a point in finite time. On the other hand,
if $V_0 \leq 0$ then the circle will monotonically shrink to a point.

\begin{table}[b!]
\center
\begin{tabular}{|r|c|c|c|c|}
\hline
$J$ 
& $\displaystyle\max_{m=0,\ldots,M} \| x(t_m) - x^m \|_{1,h}$ & EOC
& $\displaystyle\errorXtxtc$ & EOC \\ \hline
32  & 3.9796e-03 & ---  & 9.5331e-04 & ---  \\
64  & 1.0059e-03 & 1.98 & 2.4960e-04 & 1.93 \\
128 & 2.5256e-04 & 1.99 & 6.3995e-05 & 1.96 \\ 
256 & 6.3254e-05 & 2.00 & 1.6211e-05 & 1.98 \\ 
512 & 1.5827e-05 & 2.00 & 4.0803e-06 & 1.99 \\ 
1024& 3.9582e-06 & 2.00 & 1.0236e-06 & 2.00 \\ 
2048& 9.8980e-07 & 2.00 & 2.5634e-07 & 2.00 \\
\hline
\end{tabular}
\caption{Errors for the convergence test for \eqref{eq:solxpi}, 
\eqref{eq:rsol} with $r_0=1$ over the time interval $[0,1]$ 
for the scheme \eqref{eq:fdfdimpl}. 
We also display the experimental orders of convergence (EOC).
}
\label{tab:fdbimpl_T1c_p1}
\end{table}%
For the special case $V_0 = 0$, and on recalling the Gauss error function
\begin{equation*} 
\erf(z) = \tfrac2{\sqrt\pi} \int_0^z e^{-u^2} \;{\rm d}u,\quad
\erf'(z) = \tfrac2{\sqrt\pi} e^{-z^2},
\end{equation*}
we find that $r(t)$ is the solution of
\begin{equation*} 
t - \sqrt{\tfrac\pi2} r_0 \erf\left(\sqrt{\ln\frac{r_0}{r(t)}}\right) = 0,
\end{equation*}
which means that
\begin{equation} \label{eq:rsol}
r(t) = r_0 \exp(-[\erf^{-1}(\sqrt{\tfrac2\pi}\frac{t}{r_0})]^2).
\end{equation}

For the true solution
\begin{equation} \label{eq:solxpi}
 x(\rho,t) = r(t) \begin{pmatrix}
\cos g(2\pi \rho) \\ \sin g(2\pi \rho)
\end{pmatrix}, \quad g(u) = u + 0.1 \sin(u),
\end{equation}
of \eqref{eq:hcsf} we compute approximations to the errors between 
$x$ and $x^h$, the solution to \eqref{eq:fdsd} with $\mathcal{V}_0=0$, 
for the choice $r_0=1$ on the time interval $[0,1]$ with the help of the
fully discrete scheme \eqref{eq:fdfdimpl}. In particular, for the sequence of
discretization parameters $h=\frac1J = 2^{-k}$, $k=5,\ldots,11$, we let
$\Delta t = h$ and compare the grid interpolations of $x$ and $\dot x$
to their fully discrete analogues in the discrete norms \eqref{eq:discnorms}. 
These errors are reported in Table~\ref{tab:fdbimpl_T1c_p1},
where we observe the expected 
second order convergence rates from Theorem~\ref{thm:main}.

\subsection{Numerical experiments with constant initial velocity}

Throughout the remainder of the numerical results section we 
choose the discretization parameters
$J=256$ and $\Delta t = 10^{-4}$. 
Moreover, we always let $\beta=0$, unless stated otherwise.
The curve evolutions we visualize by plotting the polygonal curves 
$\Gamma^m \subset \bR^2$ defined by the vertices
$\{x^m_j\}_{j=1}^J$ and at times we also show the evolution of the length of
these curves, defined by 
\[
|\Gamma^m| = h \sum_{j=1}^J q^m_j = \sum_{j=1}^J |x^m_j - x^m_{j-1}|.
\] 
Moreover, we will often be interested in a possible blow-up in curvature, and
so we will monitor the quantity
\[
K^m_\infty = \max_{j=1,\ldots,J} \frac{|\delta \tau^m_j|}{q^m_j}
\]
as an approximation to the maximal value of 
$|\varkappa| = \frac{|\tau_\rho|}{|x_\rho|}$, recall \eqref{eq:nu}. 

In all the numerical computations in this subsection, we will
choose a constant initial velocity $\mathcal{V}_0(\rho) = V_0$.

\begin{figure}[t!]
\center
\includegraphics[angle=-90,width=0.3\textwidth]{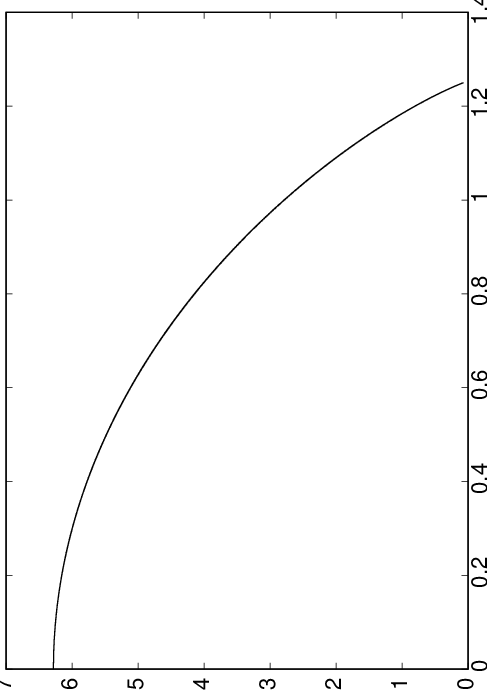}
\includegraphics[angle=-90,width=0.3\textwidth]{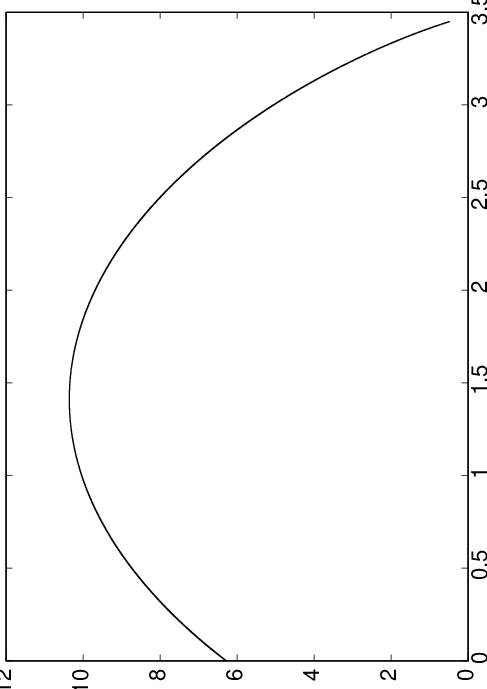}
\includegraphics[angle=-90,width=0.3\textwidth]{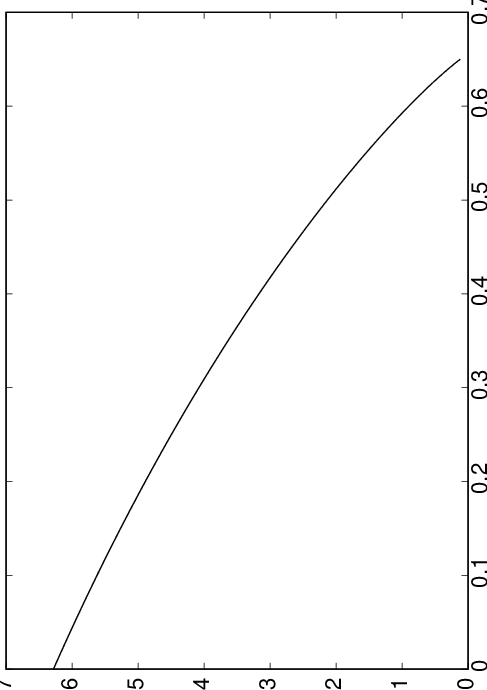}
\includegraphics[angle=-90,width=0.3\textwidth]{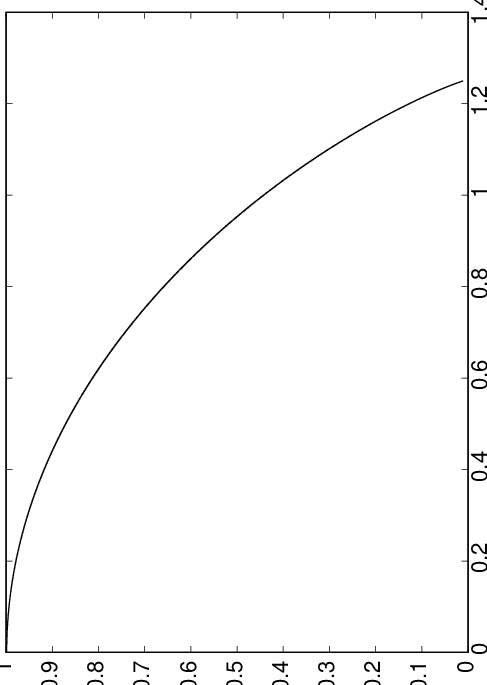}
\includegraphics[angle=-90,width=0.3\textwidth]{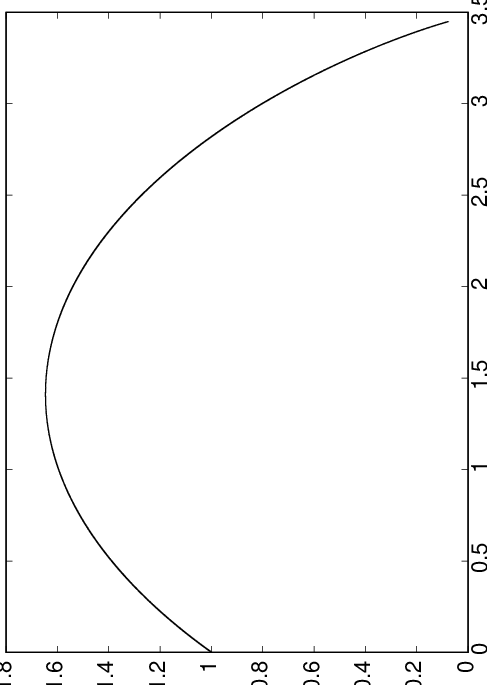}
\includegraphics[angle=-90,width=0.3\textwidth]{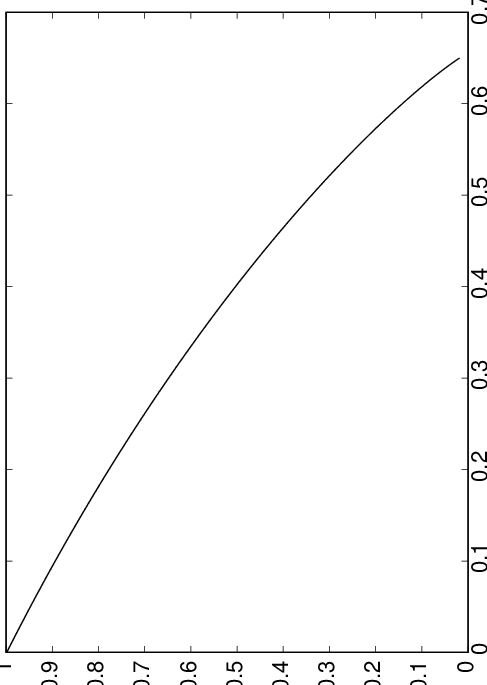}
\caption{Hyperbolic curvature flow starting from a unit circle.
Above we show the evolution of $|\Gamma^m|$ over time for 
$V_0=0$, $V_0=1$ and $V_0=-1$ (from left to right).
The final times for these computations are $T=1.25$, $T=3.45$ and $T=0.65$,
respectively.
Below we show the corresponding evolutions of $1/K^m_\infty$ over time.
}
\label{fig:circles}
\end{figure}

As discussed above, for an initial circle with uniform initial normal velocity
$\mathcal{V}_0 = V_0$, depending on the sign of $V_0 \in \bR$ the family of
circles either expands at first and then shrinks, or shrinks immediately. We
visualize these different behaviours in Figure~\ref{fig:circles}. In each case
we observe a smooth solution until the circles shrink to a point, meaning that
$|\Gamma^m|$ and $1/K^m_\infty$ approach zero at the same time.

\begin{figure}[htb]
\center
\includegraphics[angle=-90,width=0.35\textwidth]{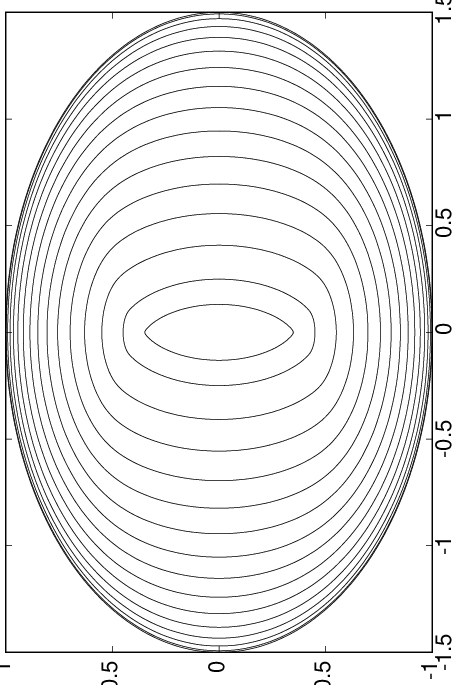}
\includegraphics[angle=-90,width=0.3\textwidth]{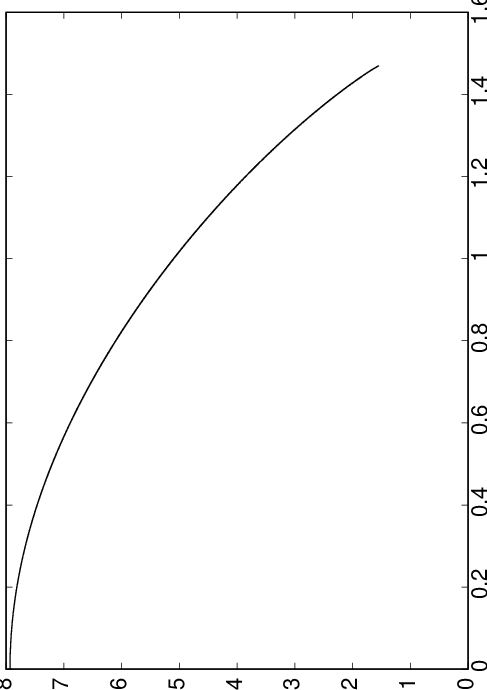}
\includegraphics[angle=-90,width=0.3\textwidth]{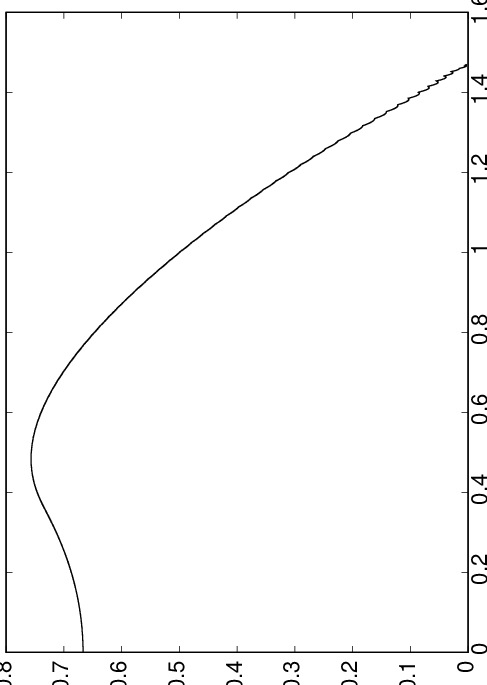}
\caption{Hyperbolic curvature flow, with $V_0=0$, starting from an
ellipse. On the left we show
$\Gamma^m$ at times $t=0,0.1,\ldots,1.4,T=1.47$.
We also show the evolutions of $|\Gamma^m|$ (middle) and 
$1/K^m_\infty$ (right) over time.
}
\label{fig:ell321V0}
\end{figure}

\begin{figure}[hbt]
\center
\includegraphics[angle=-90,width=0.35\textwidth]{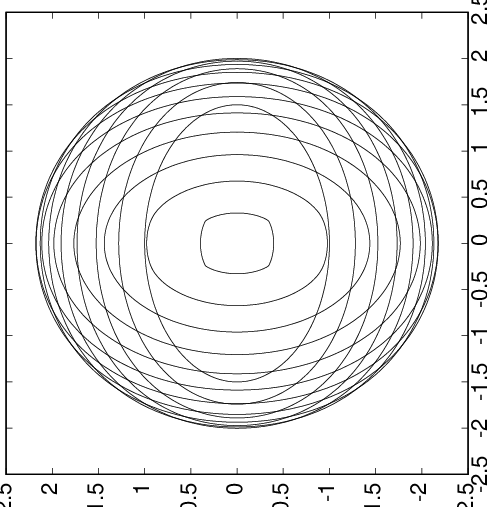}
\includegraphics[angle=-90,width=0.3\textwidth]{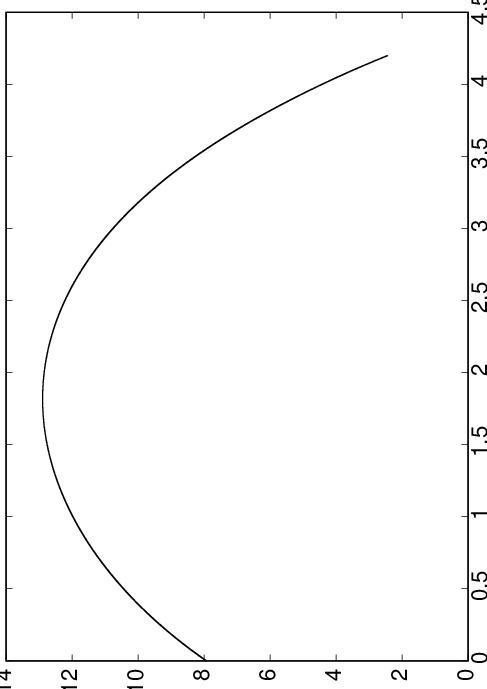}
\includegraphics[angle=-90,width=0.3\textwidth]{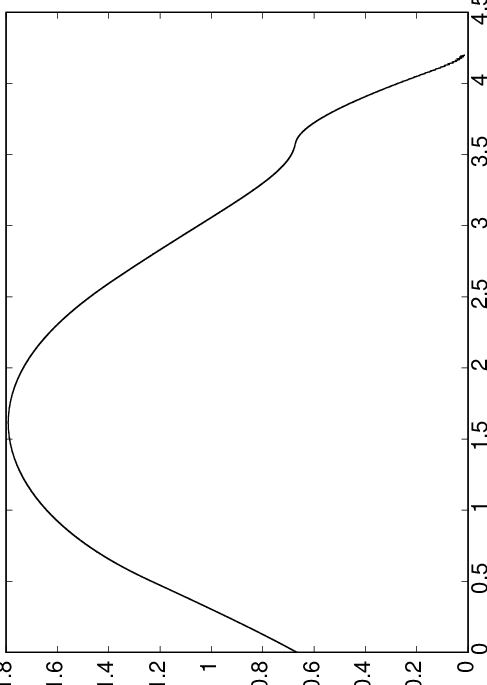}
\caption{Hyperbolic curvature flow, with $V_0=1$, starting from an
ellipse. On the left we show
$\Gamma^m$ at times $t=0,0.3,\ldots,T=4.2$.
We also show the evolutions of $|\Gamma^m|$ (middle) and 
$1/K^m_\infty$ (right) over time.
}
\label{fig:ell321V1}
\end{figure}

For the next computations we choose as initial curve a mild ellipse, with 
major axis of length 3 and minor axis of length 2.
The results for $V_0 = 0$ are shown in Figure~\ref{fig:ell321V0}, 
where we note the onset of a singularity in finite time. In particular, the
curve appears to form two kinks, leading to a blow-up in curvature.
When we choose the initial normal velocity as $V_0 = 1$, we obtain the results
shown in Figure~\ref{fig:ell321V1}. Once again we observe a blow-up in
curvature, although this time the curve does not exhibit two kinks. Instead
it seems to approach a shape with four corners.
We note that the initial ellipse at first grows towards a circle.
It then shrinks while momentarily adopting an elliptic shape, but with the
major and minor axes swapped with respect to the initial data. Towards the end
of the evolution a more circular shape appears again, which then evolves to
the limitting shape with the four corners, i.e.\ with four points where the
curvature is discontinuous.
We stress that the observed singularities in our numerical simulations 
are robust with respect to the choice of discretization parameters. For
example, refining the discretization parameters to 
$J=512$ and $\Delta t = 5 \times 10^{-5}$ gave visually indistuingishable
results compared to Figure~\ref{fig:ell321V1}.


\begin{figure}[hbt]
\center
\includegraphics[angle=-90,width=0.35\textwidth]{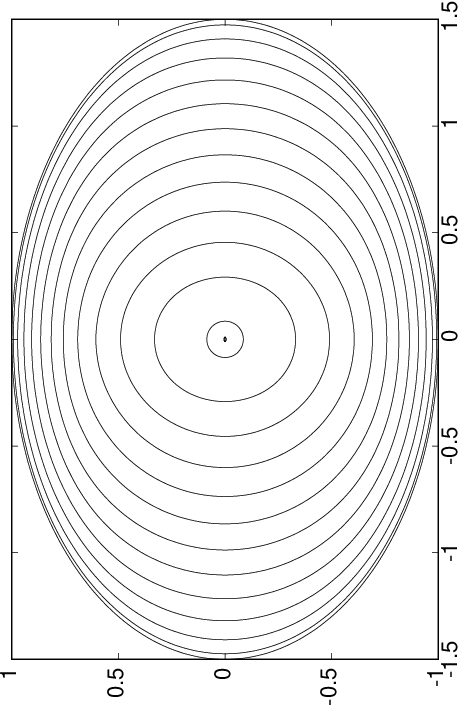}
\includegraphics[angle=-90,width=0.3\textwidth]{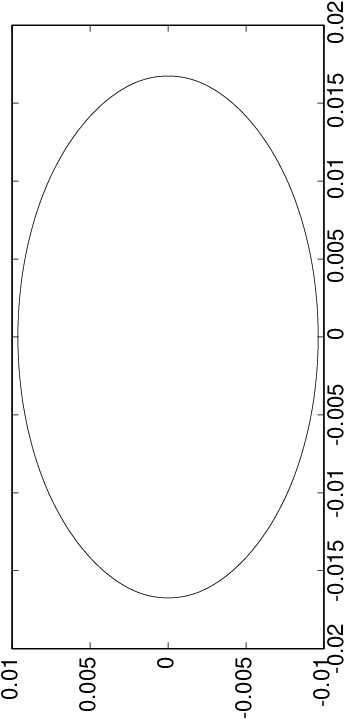}
\includegraphics[angle=-90,width=0.3\textwidth]{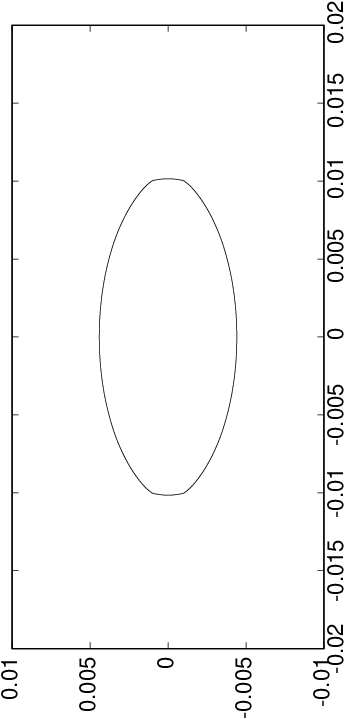}
\includegraphics[angle=-90,width=0.3\textwidth]{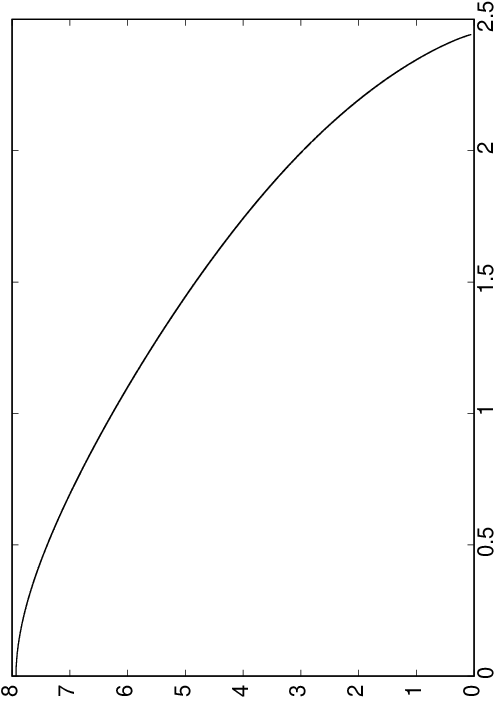}
\includegraphics[angle=-90,width=0.3\textwidth]{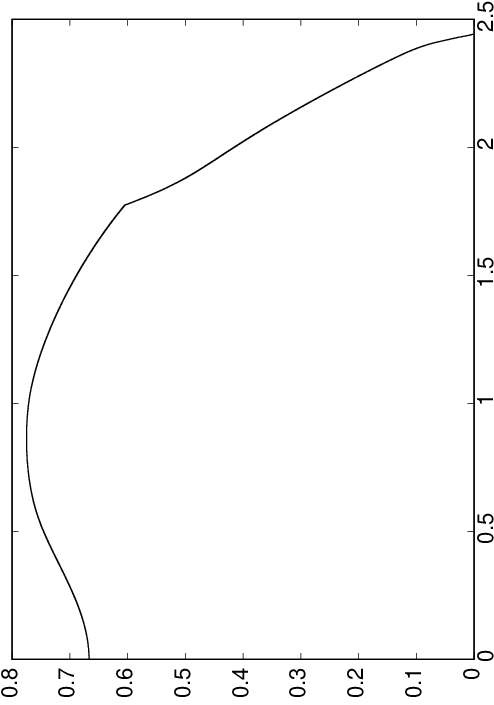}
\caption{Damped hyperbolic curvature flow, with $\beta=2$ and 
$V_0=0$, starting from an ellipse. On the top we show
$\Gamma^m$ at times $t=0,0.2,\ldots,2.4,T=2.4423$, as well as
$\Gamma^m$ separately at times $t=2.44$ and $t=T$.
Below we show the evolutions of $|\Gamma^m|$ (left) and 
$1/K^m_\infty$ (right) over time.
}
\label{fig:ell321V0b2}
\end{figure}

\begin{figure}[hbt]
\center
\includegraphics[angle=-90,width=0.35\textwidth]{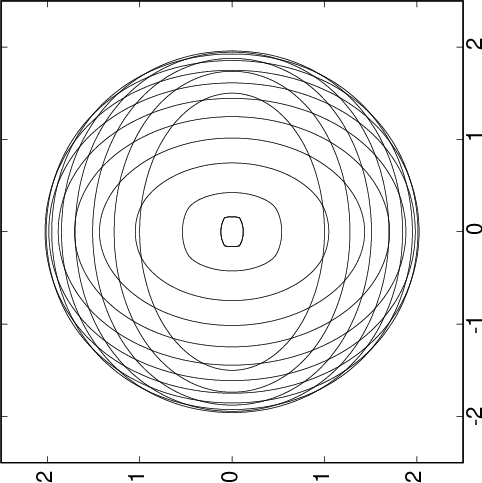}
\includegraphics[angle=-90,width=0.3\textwidth]{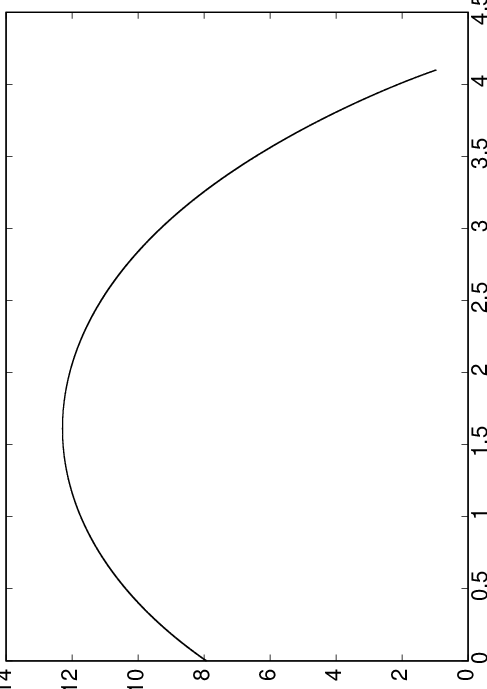}
\includegraphics[angle=-90,width=0.3\textwidth]{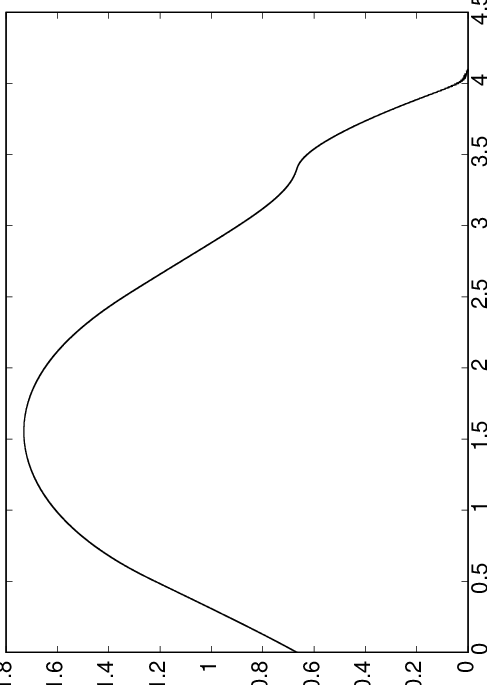}
\caption{Damped hyperbolic curvature flow, with $\beta=0.1$ and 
$V_0=1$, starting from an ellipse. On the left we show
$\Gamma^m$ at times $t=0,0.3,\ldots,3.9,T=4.1$.
We also show the evolutions of $|\Gamma^m|$ (middle) and 
$1/K^m_\infty$ (right) over time.
}
\label{fig:ell321V1b01}
\end{figure}

We are interested in the effect of the parameter $\beta$ on these developing
singularities, and would expect some damping or smoothing to be observable
for $\beta > 0$.
Repeating the simulation from Figure~\ref{fig:ell321V0}
with $\beta=2$ yields the results in Figure~\ref{fig:ell321V0b2},
where we observe that the blow-up in curvature now happens much later, 
when the curve itself is almost extinct. We also see a marked change in the
profile of the evolving curve. While in Figure~\ref{fig:ell321V0} at late
times the curve resembles an ellipsoid aligned with the $x_2$-axis, 
the evolution in Figure~\ref{fig:ell321V0b2} for long times appears to 
approach a circle, until towards the very end it starts to resemble an 
ellipsoid aligned with the $x_1$-axis.
In addition, a repeat of Figure~\ref{fig:ell321V1} now with $\beta=0.1$ is 
shown in Figure~\ref{fig:ell321V1b01}, where once again we note that 
visually the curve appears smoother for longer, 
until eventually the curvature blows up due to facetting on the left and right
sides of the curve.

\begin{figure}[hbt]
\center
\includegraphics[angle=-90,width=0.6\textwidth]{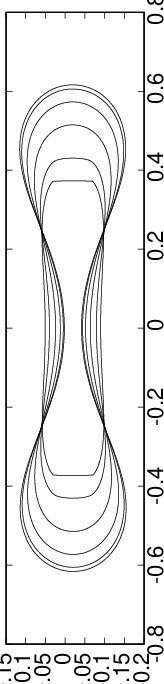}\\
\includegraphics[angle=-90,width=0.3\textwidth]{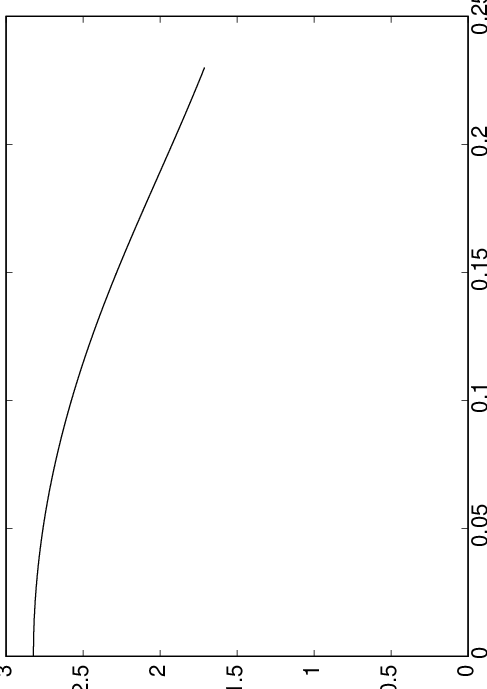}
\includegraphics[angle=-90,width=0.3\textwidth]{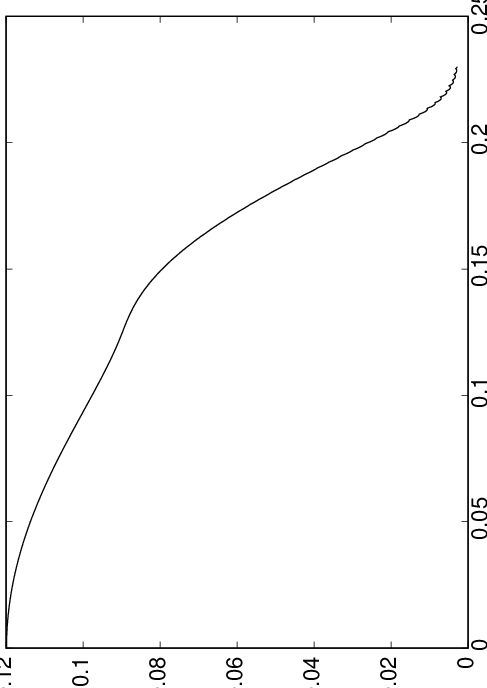}
\caption{Hyperbolic curvature flow, with $V_0=0$, starting from a
smooth dumbbell. 
On top we show $\Gamma^m$ at times $t=0,0.05,\ldots,0.2,T=0.23$.
Below we show the evolutions of $|\Gamma^m|$ (left) and 
$1/K^m_\infty$ (right) over time.
}
\label{fig:dumbbellV0}
\end{figure}

\begin{figure}[hbt]
\center
\includegraphics[angle=-90,width=0.3\textwidth]{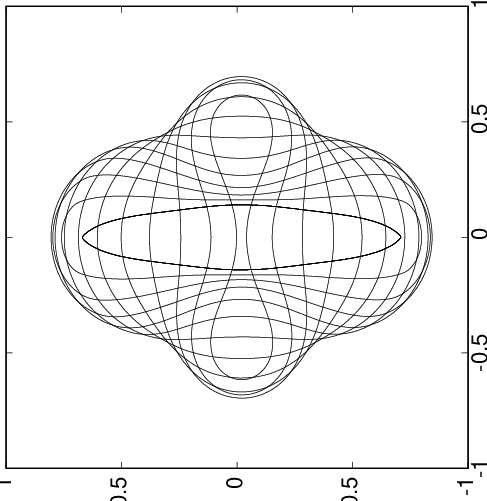}
\includegraphics[angle=-90,width=0.3\textwidth]{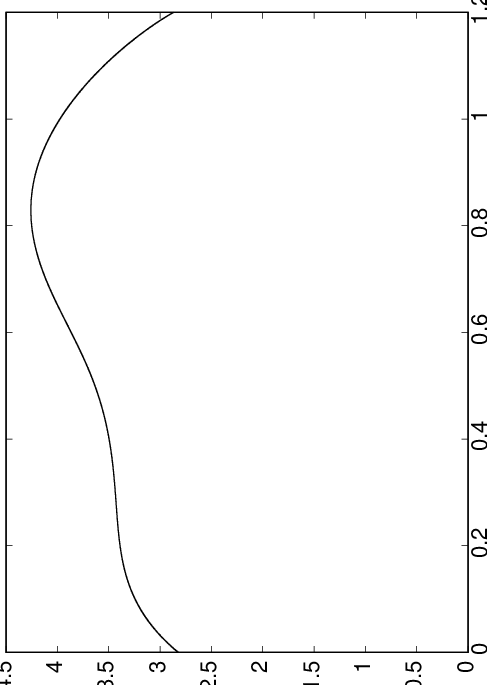}
\includegraphics[angle=-90,width=0.3\textwidth]{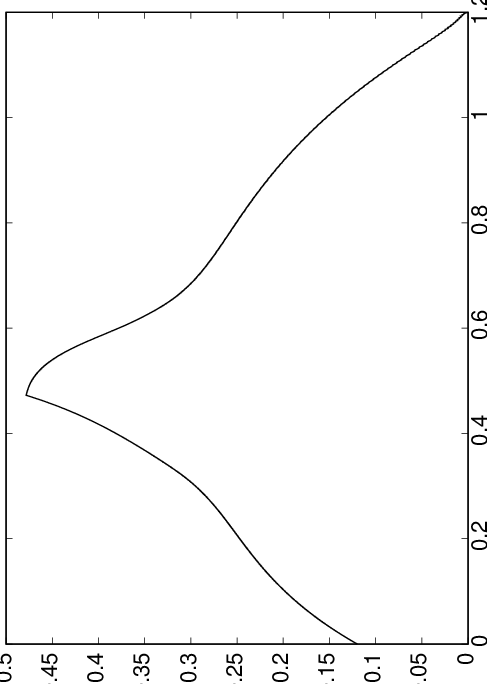}
\mbox{
\includegraphics[angle=-90,width=0.22\textwidth]{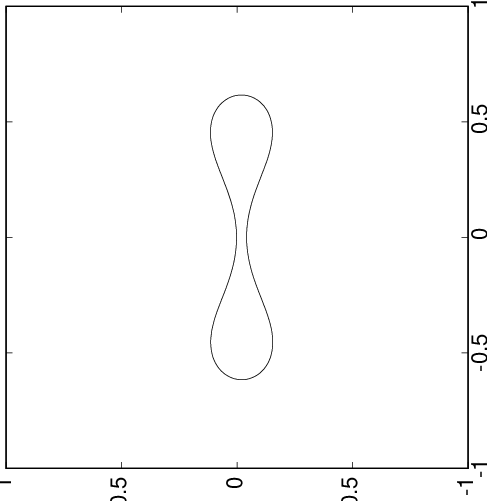}
\includegraphics[angle=-90,width=0.22\textwidth]{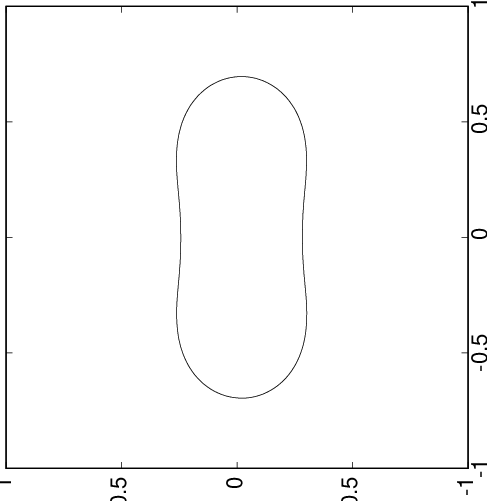}
\includegraphics[angle=-90,width=0.22\textwidth]{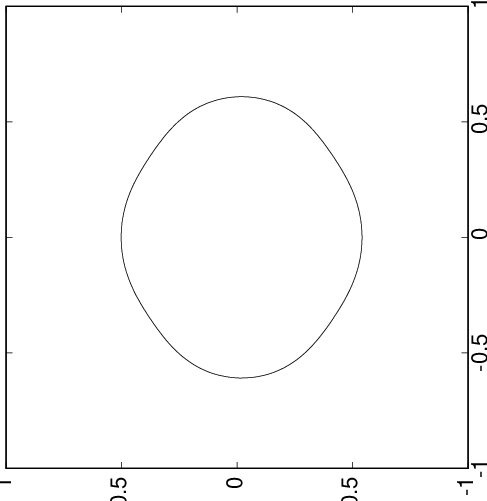}
\includegraphics[angle=-90,width=0.22\textwidth]{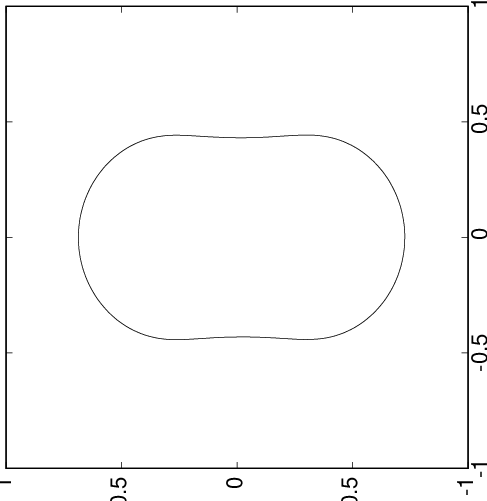}
}
\mbox{
\includegraphics[angle=-90,width=0.22\textwidth]{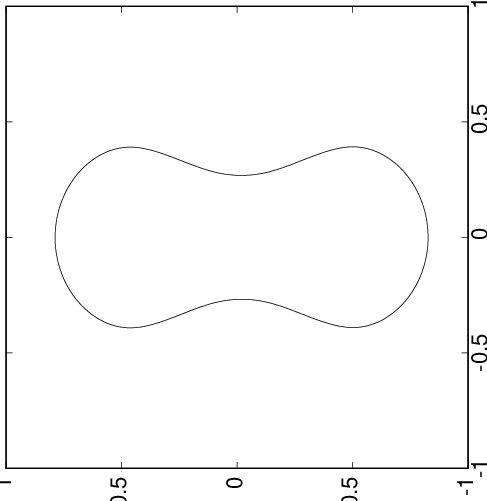}
\includegraphics[angle=-90,width=0.22\textwidth]{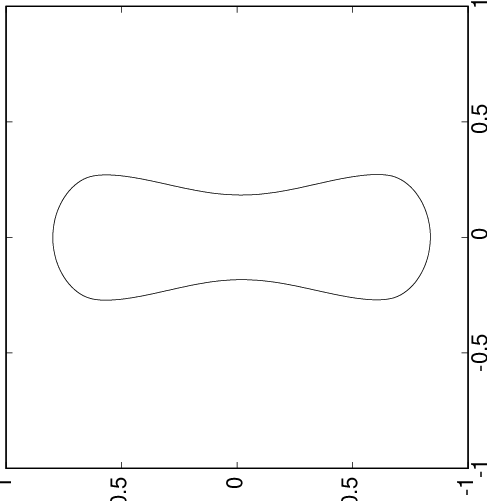}
\includegraphics[angle=-90,width=0.22\textwidth]{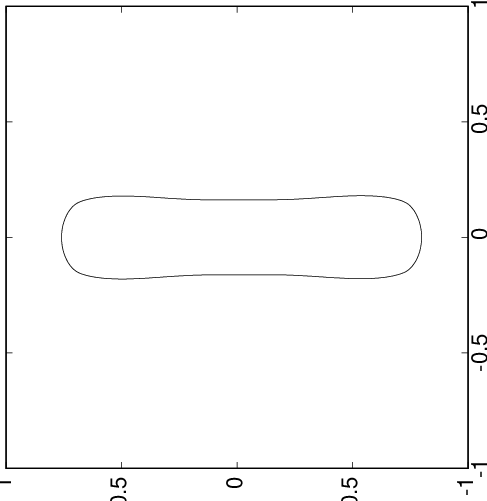}
\includegraphics[angle=-90,width=0.22\textwidth]{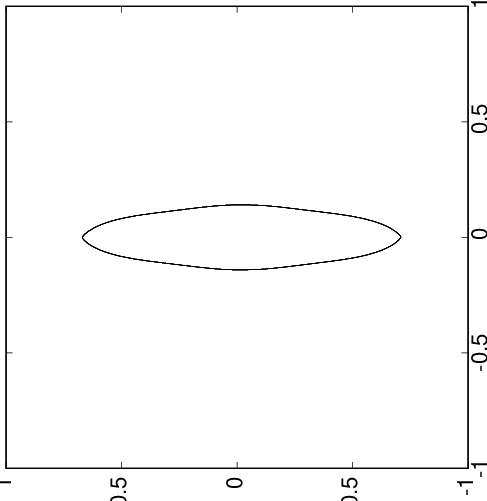}
}
\caption{Hyperbolic curvature flow, with $V_0=1$, starting from a
smooth dumbbell. 
On top we show $\Gamma^m$ at times $t=0,0.1,\ldots,T=1.2$ (left), as well as
evolutions of $|\Gamma^m|$ (middle) and $1/K^m_\infty$ (right) over time. 
Below we visualize $\Gamma^m$ separately
at times $t=0,0.2,0.4,0.6$ and $0.8,1,1.1,1.2$.
}
\label{fig:dumbbellV1}
\end{figure}

\begin{figure}[hbt]
\center
\includegraphics[angle=-90,width=0.6\textwidth]{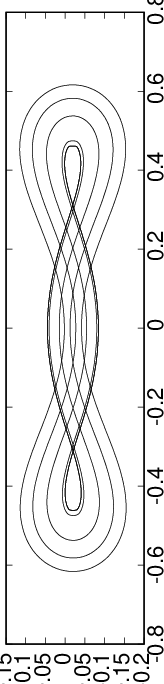}\\
\includegraphics[angle=-90,width=0.3\textwidth]{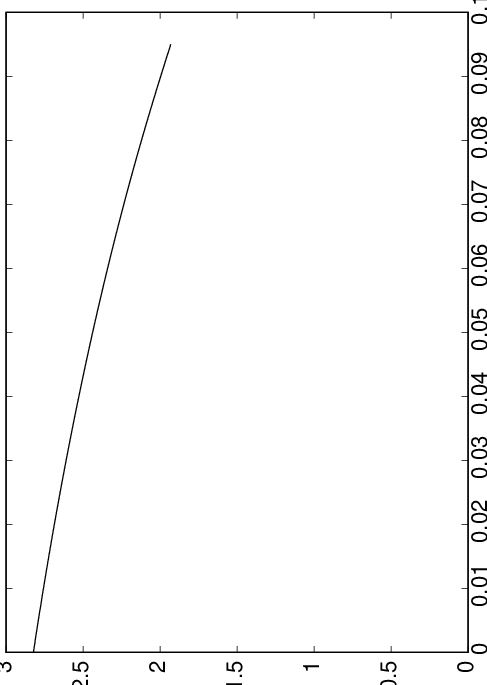}
\includegraphics[angle=-90,width=0.3\textwidth]{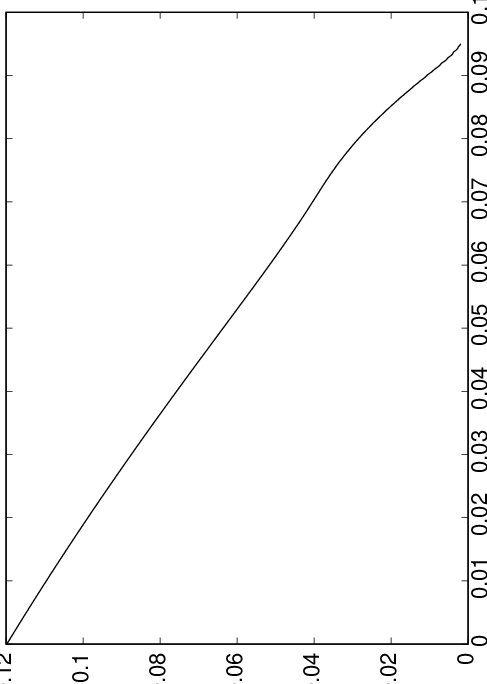}
\caption{Hyperbolic curvature flow, with $V_0=-1$, starting from a
smooth dumbbell. 
On top we show $\Gamma^m$ at times $t=0,0.03,\ldots,0.09,T=0.095$.
Below we show the evolutions of $|\Gamma^m|$ (left) and 
$1/K^m_\infty$ (right) over time.
}
\label{fig:dumbbellV-1}
\end{figure}

Finally, we also consider some numerical experiments where the initial data is
nonconvex. For the simulation in Figure~\ref{fig:dumbbellV0} we start from a
smooth dumbbell-like initial curve. We observe that the curve starts to shrink
until it eventually exhibits two facets on the left and right, which leads to a
blow-up in the curvature. Repeating the simulation for the constant 
initial velocity $V_0=1$ yields the results in Figure~\ref{fig:dumbbellV1}.
Now the curve first expands vertically into a convex curve that expands 
further, until it narrows on the $x_1$-axis towards the origin to create a 
new nonconvex shape that resembles a variant of the initial data that is 
now aligned with the $x_2$-axis. At this stage the curve begins again to 
expand into a convex shape that then shrinks until two developing kinks 
at the top and bottom of the curve lead to a blow-up in curvature.
Interestingly, when we use the initial velocity $V_0=-1$ the curve soon
self-intersects, see Figure~\ref{fig:dumbbellV-1},
which is something the parametric formulation is blind towards. 
Similarly to the evolution in Figure~\ref{fig:dumbbellV0}, the solution
approaches a blow-up in curvature when two facets are about to be created on
the left and right sides of the curve.

In conclusion we remark that the onset of a blow-up in curvature in finite
time for strictly convex initial data as observed in Figure~\ref{fig:ell321V0}
confirms the theoretical predictions in \cite{KongLW09}. In addition,
Figure~\ref{fig:ell321V1} demonstrates that the same can be observed for
an outward initial velocity $\mathcal{V}_0 \nu(\cdot,0)$.
Finally, from our remaining numerical simulations we
conjecture that also nonconvex initial data can exhibit the same phenomenon.

\subsection{Numerical experiments with nonconstant initial velocity}
\label{sec:nr54}

\begin{figure}[hbt]
\center
\includegraphics[angle=-90,width=0.35\textwidth]{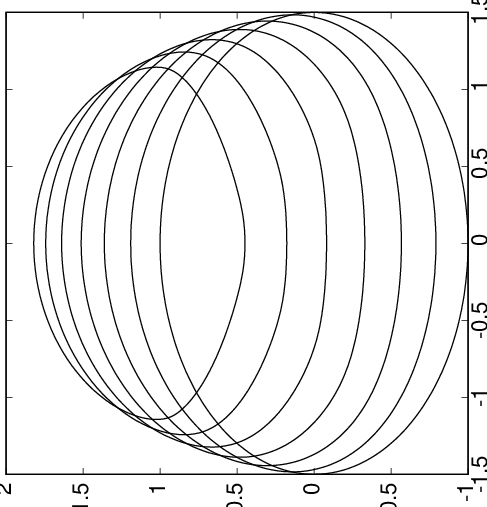}
\includegraphics[angle=-90,width=0.3\textwidth]{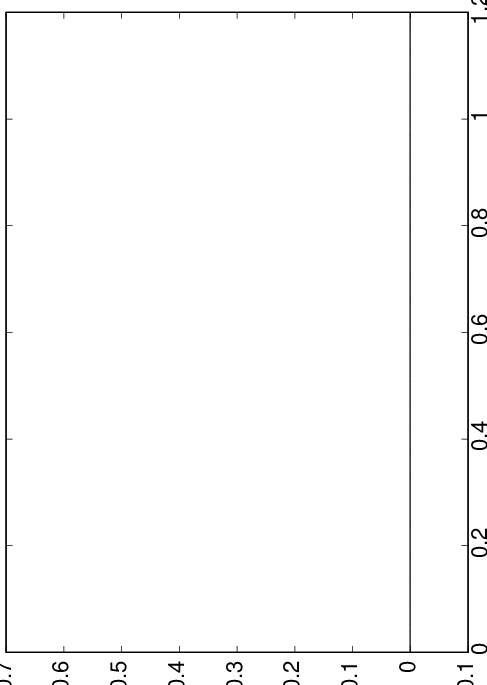}
\includegraphics[angle=-90,width=0.3\textwidth]{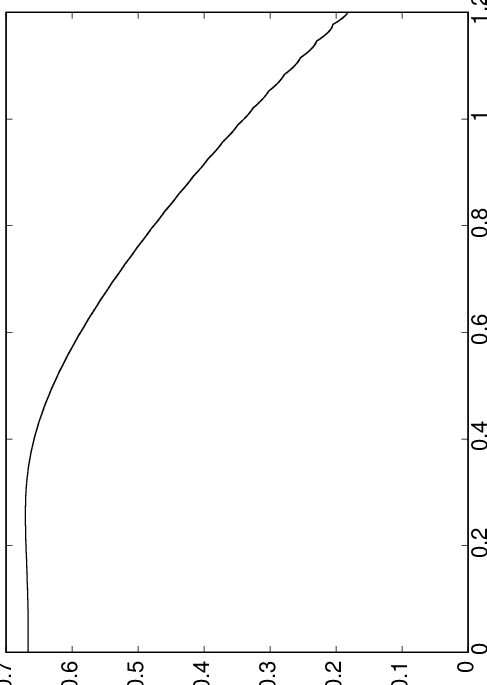}
\caption{Hyperbolic curvature flow, with $\mathcal{V}_0(\rho)=
\sin(2 \pi \rho)$, starting from an ellipse parameterized by
$x_0(\rho) = (\frac32\cos(2 \pi \rho) , \sin(2 \pi \rho))^T$. 
On the left we show $\Gamma^m$ at times $t=0,0.2,\ldots,T=1.2$.
We also show the evolutions of $\| D_t x^{m+1} \cdot \theta^m \|_{0,h}$ 
(middle) and $1/K^m_\infty$ (right) over time.
}
\label{fig:ell321specialV2}
\end{figure}

In this final subsection we report on a numerical simulation
with a nonconstant initial velocity $\mathcal{V}_0$. In particular,
we repeat the experiment from Figure~\ref{fig:ell321V1}, but now choose
$\mathcal{V}_0(\rho) = \sin(2 \pi \rho)$, with $x_0(\rho) = 
(\frac32\cos(2 \pi \rho) , \sin(2 \pi \rho))^T$. The evolution can be seen
in Figure~\ref{fig:ell321specialV2}. Note that due to the given
initial velocity, the curve rises and shrinks at the same time. Towards the end
of the evolution a flat patch appears to develop at the bottom part of the
curve. For a later comparison, we also provide a plot of the discrete
tangential velocity
\[
\| D_t x^{m+1} \cdot \theta^m \|_{0,h} 
: = \left( h \sum_{j=1}^J \left| \frac{x^{m+1}_j - x^m_j}{\Delta t} \cdot \theta^m_j \right|^2 \right)^{\frac{1}{2}}.
\]
over time in Figure~\ref{fig:ell321specialV2}. Since \eqref{eq:fdfdimpl} is a
discrete approximation of the normal flow \eqref{eq:hcsf}, the quantity stays
nearly equal to zero throughout the evolution. 

\begin{figure}[hbt]
\center
\includegraphics[angle=-90,width=0.35\textwidth]{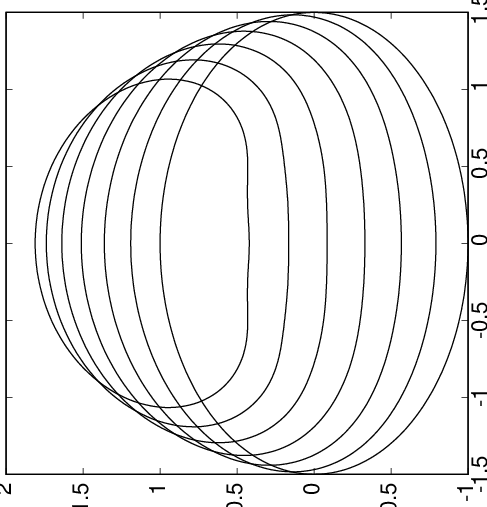}
\includegraphics[angle=-90,width=0.3\textwidth]{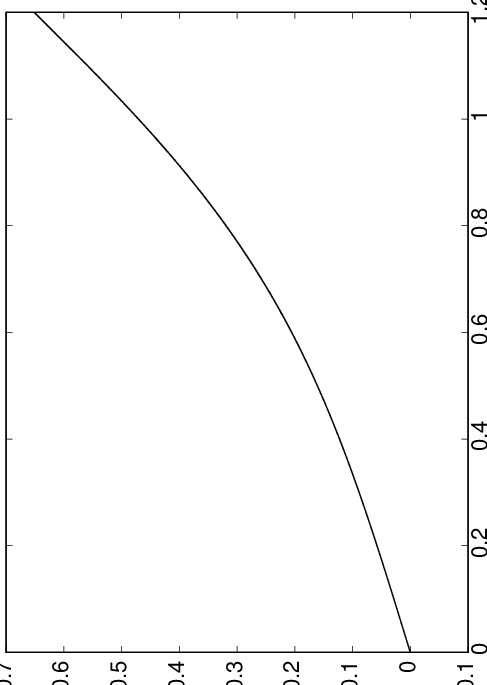}
\includegraphics[angle=-90,width=0.3\textwidth]{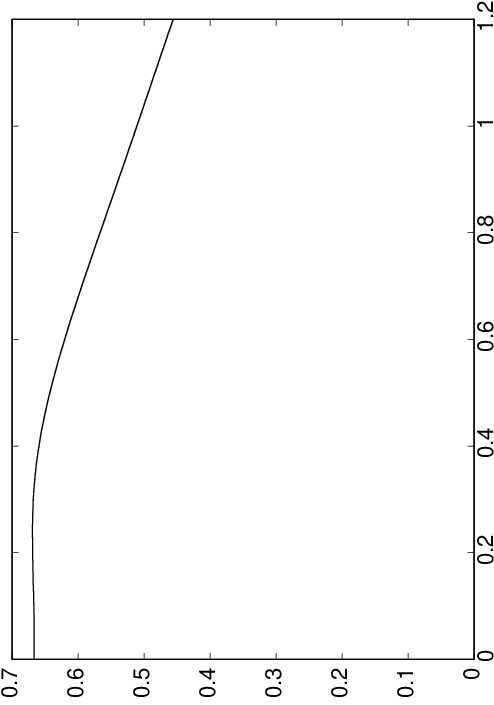}
\caption{The flow \eqref{eq:xttxss0} with \eqref{eq:init} 
for $\mathcal{V}_0(\rho)=\sin(2 \pi \rho)$, starting from an ellipse 
parameterized by $x_0(\rho) = (\frac32\cos(2 \pi \rho) , \sin(2 \pi \rho))^T$. 
On the left we show $\Gamma^m$ at times $t=0,0.2,\ldots,T=1.2$.
We also show the evolutions of $\| D_t x^{m+1} \cdot \theta^m \|_{0,h}$ 
(middle) and $1/K^m_\infty$ (right) over time.
}
\label{fig:debugell321specialV2}
\end{figure}

We mentioned in the introduction that a question of mathematical interest is
whether solutions to \eqref{eq:xttxss0} parameterize curves evolving according
to \eqref{eq:Gurtinbeta0}. We now provide some numerical evidence that this is
not the case. In order to numerically approximate solutions to
\eqref{eq:xttxss0}, we naturally adapt the scheme \eqref{eq:fdfdimpl},
for $\beta = 0$, by omitting the last term on the right hand side of
\eqref{eq:fdfdimpl}. For this new scheme we then repeat the computation from
Figure~\ref{fig:ell321specialV2} using exactly the same discrete initial data.
The ensuing evolution, shown in Figure~\ref{fig:debugell321specialV2},
is close to what we observed before, but ultimately
differs. The differences are most pronounced in the final shape of $\Gamma^m$
and in the plot of $1/K^m_\infty$ over time.
We remark that a main difference between \eqref{eq:hcsf} and \eqref{eq:xttxss0} 
is that the former is a normal flow, while the latter allows for a nonzero
tangential component of the velocity $x_t$. Once again this is confirmed by our
numerical experiment, as can be seen from the plot of
$\| D_t x^{m+1} \cdot \theta^m \|_{0,h}$ in 
Figure~\ref{fig:debugell321specialV2}, which seems to be monotonically
increasing.
We remark that we repeated the simulations in Figures~\ref{fig:ell321specialV2}
and \ref{fig:debugell321specialV2} with finer discretization parameters and
obtained visually indistinguishable results. Hence we are confident that
the displayed evolution provide numerical evidence that the two PDEs
\eqref{eq:abxttxss} and \eqref{eq:xttxss0}, with the initial conditions
\eqref{eq:init}, parameterize different curve evolutions.

\def\soft#1{\leavevmode\setbox0=\hbox{h}\dimen7=\ht0\advance \dimen7
  by-1ex\relax\if t#1\relax\rlap{\raise.6\dimen7
  \hbox{\kern.3ex\char'47}}#1\relax\else\if T#1\relax
  \rlap{\raise.5\dimen7\hbox{\kern1.3ex\char'47}}#1\relax \else\if
  d#1\relax\rlap{\raise.5\dimen7\hbox{\kern.9ex \char'47}}#1\relax\else\if
  D#1\relax\rlap{\raise.5\dimen7 \hbox{\kern1.4ex\char'47}}#1\relax\else\if
  l#1\relax \rlap{\raise.5\dimen7\hbox{\kern.4ex\char'47}}#1\relax \else\if
  L#1\relax\rlap{\raise.5\dimen7\hbox{\kern.7ex
  \char'47}}#1\relax\else\message{accent \string\soft \space #1 not
  defined!}#1\relax\fi\fi\fi\fi\fi\fi}

\end{document}